\documentclass[11pt,a5paper]{article}
\usepackage[a5paper,margin=10mm,tmargin=15mm]{geometry}

\title{Adaptively detect and accurately resolve macro-scale shocks in an efficient Equation-Free multiscale simulation}

\author{John~Maclean\thanks{\protect\raggedright%
School of Mathematical Sciences, University of Adelaide, South Australia.
\url{http://www.adelaide.edu.au/directory/john.maclean}}
\and J.~E.~Bunder\thanks{\protect\raggedright%
School of Mathematical Sciences, University of Adelaide, South Australia.
\protect\url{mailto:judith.bunder@adelaide.edu.au},
\protect\url{http://orcid.org/0000-0001-5355-2288}}
\and 
I.~G. Kevrekidis
\thanks{Departments of Chemical and Biomolecular Engineering and Applied Mathematics and Statistics, Johns Hopkins University, Baltimore, Maryland, USA.
\protect\url{https://orcid.org/0000-0003-2220-3522}}
\and A.~J.~Roberts\thanks{\protect\raggedright%
School of Mathematical Sciences, University of Adelaide, South Australia.
\url{https://profajroberts.github.io},
\url{http://orcid.org/0000-0001-8930-1552}}
}

\date{\today}

\usepackage[T1]{fontenc}
\usepackage[dvipsnames]{xcolor}
\usepackage{enumitem}

\usepackage{pgfplots,ifthen} 
\pgfplotsset{compat=newest} 
\usepgfplotslibrary{external} 
\usetikzlibrary{decorations.markings}
\usetikzlibrary{shapes,arrows,fit}
\usetikzlibrary{positioning}
\usetikzlibrary{calc} 

\usepackage{natbib}
\AtEndDocument{%
    
    \IfFileExists{ajr.sty}%
        {\bibliography{bibexport,ajr,bib}}%
        {\bibliography{bibexport}}%
    }

\usepackage{url,microtype,amsmath,amssymb,graphicx,hyperref,doi,defns}
\hypersetup{colorlinks,allcolors=RoyalBlue,breaklinks,%
            bookmarks,bookmarksopen}
\makeatletter
\AtBeginDocument{{\def\and{, }\def\thanks#1{}%
  \hypersetup{
    pdfauthor={\@author},
    pdftitle={\@title}}}
    }
\makeatother
\def\Pat{\RaisedName{Pat}meso-patch} 
\def\sige{\RaisedName{sige}\sigma_\epsilon}
\def\sigg{\RaisedName{sigg}\sigma_\gamma}
\def\dx{\RaisedName{dx}d}
\Vec u  
\Bb R \Bb N
\Cal N
\RaisedNamesfalse

\def\mrg{'\RaisedName{mrg}}

\usepackage{mycleveref}

\graphicspath{{figs/}}

\definecolor{Mcol}{RGB}{0,  73,   183}
\definecolor{Mtxt}{RGB}{0,  65,   165} 
\definecolor{mcol}{RGB}{217,  95,   2} 
\definecolor{mtxt}{RGB}{200,  80,   2} 
\definecolor{excol}{RGB}{76,  0,   153}  
\definecolor{extxt}{RGB}{80,  0,   163}  

\begin{document}

\maketitle

\begin{abstract}
The Equation-Free approach to efficient multiscale numerical computation marries trusted micro-scale simulations to a framework for numerical macro-scale reduction---the patch dynamics scheme.
A recent novel patch scheme empowered the Equation-Free approach to simulate systems containing shocks on the macro-scale.
However, the scheme did not predict the formation of shocks accurately, and it could not simulate moving shocks. 
This article resolves both issues, as a first step in one spatial dimension, by embedding the Equation-Free, shock-resolving patch scheme within a classic framework for adaptive moving meshes. 
Our canonical micro-scale problems exhibit heterogeneous nonlinear advection and heterogeneous diffusion.
We demonstrate many remarkable benefits from the moving patch scheme, including efficient and accurate macro-scale prediction despite the unknown macro-scale closure.
Equation-free methods are here extended to simulate moving, forming and merging shocks without a priori knowledge of the existence or closure of the shocks.
Whereas adaptive moving mesh equations are typically stiff, typically requiring small time-steps on the macro-scale, the moving macro-scale mesh of patches is typically not stiff given the context of the micro-scale time-steps required for the sub-patch dynamics.
\end{abstract}


\section{Introduction}

Engineering and science is increasingly operating at the micro- and nano-scales to achieve objectives at the macro-scale \cite[]{Ashby2010}.
As commented in the \textsc{siam} review by \cite{Rude2018}, ``in science and engineering simulations, large differences in temporal and spatial scales must be resolved \ldots\ Challenging features common in applications include inhomogeneity [heterogeneity]''. 
For example, the dynamic nature of crack propagation and the resulting discontinuities in material structure have created many mathematical approaches \cite[e.g.]{Tadmor1996, Guenter2011}.
We here develop novel and efficient multiscale computational methods based upon the so-called \emph{Equation-Free patch scheme} \cite[e.g.]{Kevrekidis03b, Samaey04, Samaey08}, for the sparse simulation of multiscale spatial systems exhibiting both heterogeneity and a finite number of localised cracks\slash fronts\slash shocks. 
The Equation-Free methodology has remarkably efficient scaling for multiscale systems in multiple spatial dimensions. 
The computational speedup for a multiscale system in $n$-D space, with scale separation parameter~$\varepsilon$, is~\Ord{\varepsilon^{-n}} as discussed for efficient algorithms implementing Equation-Free algorithms in 2D and 3D~space \cite[]{Roberts2011a, Maclean2020c}, and including 2D wave systems \cite[\S2.2.2]{Bunder21}.

\begin{figure}
\centering
\begin{tabular}{@{}c c@{}}
(a) & (b) \\
\includegraphics[width=0.48\linewidth]{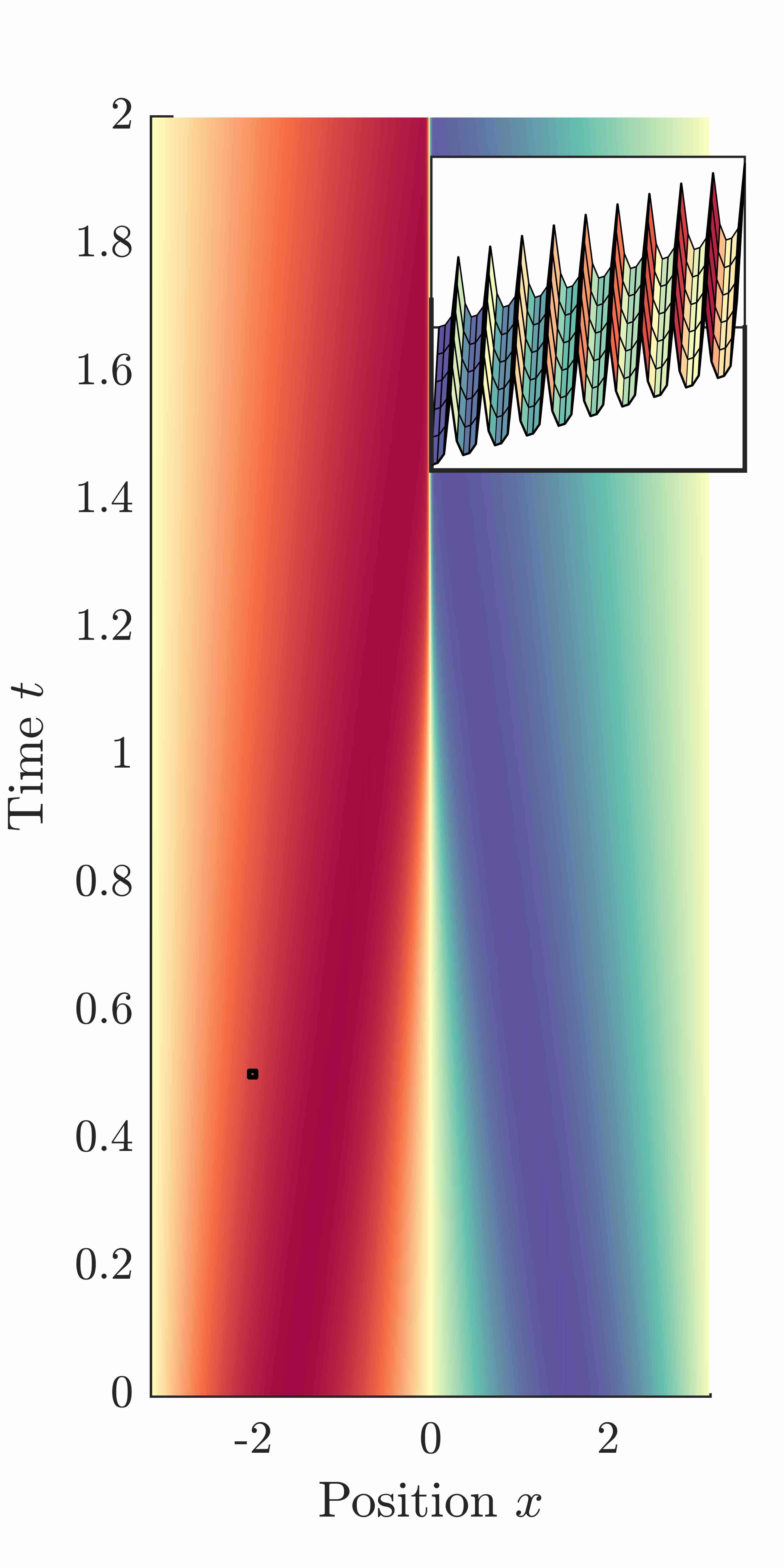} &
\includegraphics[width=0.48\linewidth]{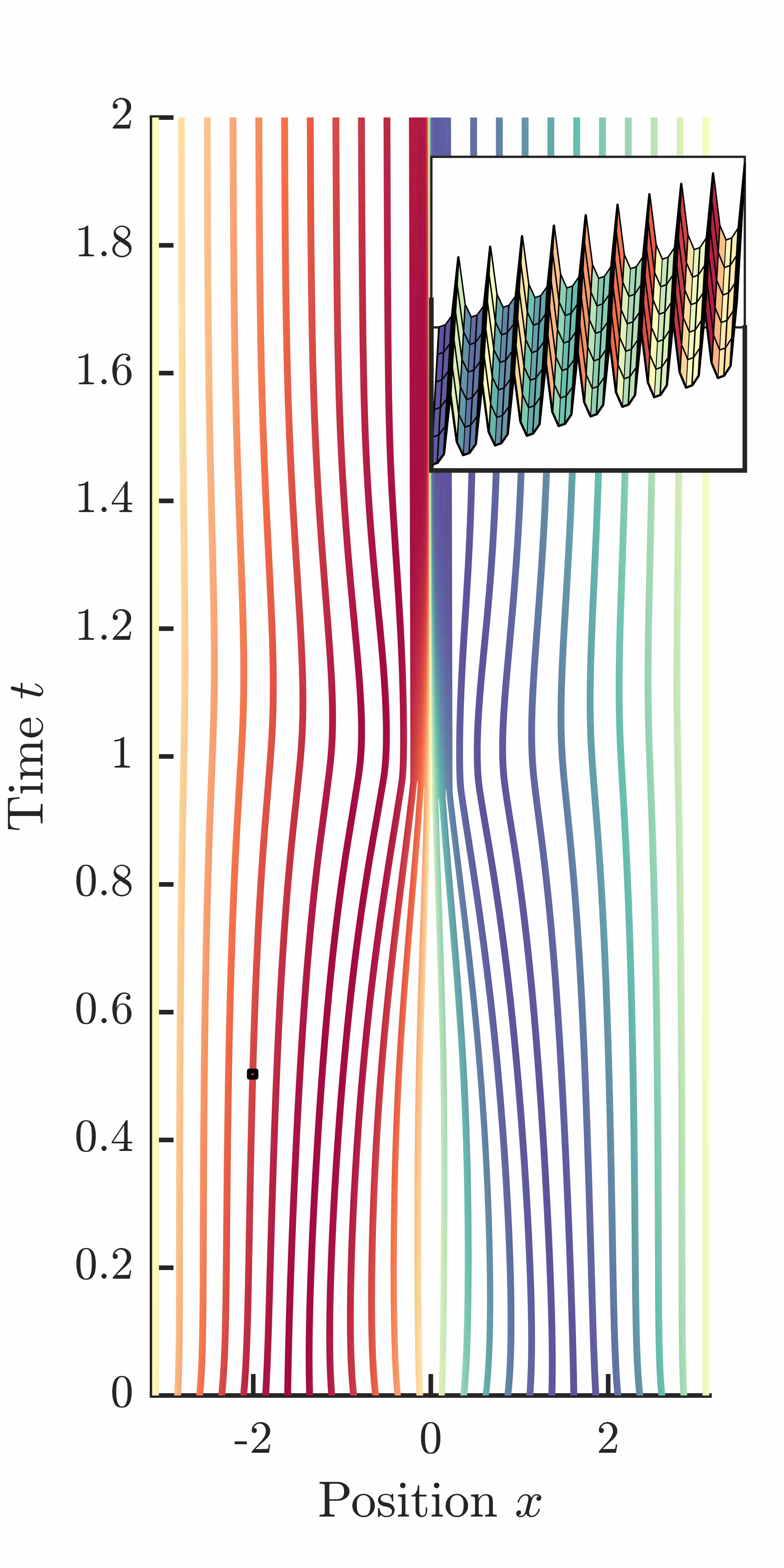}
\end{tabular}       
\caption{Simulations resolve an emergent shock\slash front in the heterogeneous micro-scale system~\cref{hetModel}.
Shading and colour show a local average of the computed field~$u(x,t$).
\quad(a)~An integration of~\cref{hetModel} over the entire domain is computationally expensive.
The inset at top shows a small part of the field, that of the black box located near $(x,t)=(-2,0.5)$, reveals the five-periodic micro-scale structure.
\quad(b)~Simulation by moving patches that cover just one quarter of the domain.
As the shock forms, six patches merge to form a so-called \emph{\Pat} (\cref{sec:movDpat}), empowering accurate resolution of the shock.
The inset confirms that the scheme accurately reproduces the \text{heterogeneous micro-scale.}}
    \label{fig:hetsin}
\end{figure}

Using established moving mesh methods \cite[e.g.]{Li1998, MacKenzie07, Budd2009, Huang10}, this article starts developing novel moving and merging small spatial patches in order to efficiently resolve both the shocks and the otherwise smooth macroscale. 
This novel moving patch scheme should readily generalise to multi-D space analogous to established moving meshes.
This article takes the first step by developing the techniques for just 1D~space in order to create a foundation to subsequently developing techniques for multi-D space.
Hence we here consider the evolution of fields~\(u(x,t)\) in 1D space~\(x\) and time~\(t\).%
\footnote{The term ``space'' need not necessarily mean physical space, just the domain over which the solution is defined at any time instant. For example, the solution could be a probability distribution, or in a frequency domain, or a velocity-space domain.}
But for computation we restrict attention to values~\(u_k(t)\) on the micro-scale grid~\(x_k\) with micro-scale equi-spacing~\(\dx\) (e.g., a lattice Boltzmann or finite difference\slash element\slash volume method).
Inspired by the heterogeneous Burgers' \pde, \(u_t=\partial_x\big[\epsilon(x)u_x-\gamma(x)u^2\big]\), we focus on the following prototype system possessing micro-scale heterogeneous advection and heterogeneous diffusion:
\begin{align}
    \frac{du_k}{dt} &= \frac1{\dx^2}\left[ \epsilon_k(u_{k+1}-u_k) -\epsilon_{k-1}(u_k-u_{k-1}) \right] 
\nonumber\\&\quad{}
 - \frac1{2\dx}\left[ \gamma_{k+1} u_{k+1}^2 -  \gamma_{k-1} u_{k-1}^2 \right] .
  \label{hetModel}
\end{align}
We assume herein that the heterogeneities~$\epsilon_k,\gamma_k$ are periodic in~\(k\), say \(\kappa\)-periodic in~\(k\) (\(\kappa=1\) if homogeneous), so that the system is \(\kappa\dx\)-periodic in~\(x\).
\Cref{fig:hetsin}(a) shows an example simulation with macro-scale modulation of micro-scale variations in the field~$u_k$, and the formation of a macro-scale shock\slash front in the field (at \(x=0\)). 
The scenario is that there no shocks to the micro-scale system, the shocks only appear when solutions are viewed on a macro-scale: as discussed in \cref{sec:e1}, the emergent macro-scale shocks we discuss are resolved accurately over 50~or so micro-scale grid points in space.
The prototype system of \ode{}s~\cref{hetModel} is significant as a generic member of its universality class of multiscale, heterogeneous, nonlinear, shock-forming, problems in~1D for which we have a trustworthy \text{micro-scale computational system.}


The importance of the prototype heterogeneous system~\cref{hetModel} is the following.
Whereas it is plausible that~\cref{hetModel} possess an effective homogenised model on a macro-scale, there is no known homogenised model.
Further, there is no known model of a Rankine--Hugoniot like condition for any emergent macro-scale shocks in~\cref{hetModel}.
To circumvent the absence of such macro-scale models, we here develop extensions to the multiscale Equation-Free patch scheme \cite[e.g.]{Kevrekidis03b, Samaey04, Samaey08} to efficiently and accurately simulate nonlinear heterogeneous systems with shocks, such as~\cref{hetModel}---developments which we contend will apply to general systems with complicated multi-physics \text{heterogeneous micro-scales.}

The fundamental idea of the patch scheme is that small patches of a given complicated micro-scale simulation (e.g.,~\cref{hetModel}) are well separated in space (e.g., \cref{fig:hetsin}(b)), and linked via mathematically designed coupling \cite[]{Roberts00a}, to efficiently cross the scales in order to accurately predict variations in space-time of macro-scale quantities of interest \cite[e.g.]{Roberts2014a, Bunder2013c}.
Other cognate numerical multiscale techniques include the molecular dynamics of \cite{car85}, numerical\slash computational homogenization \cite[e.g.,][]{Saeb2016, Geers2017, Peterseim2019}, Heterogeneous Multiscale Methods \cite[reviewed by][]{e_heterogeneous_2007, abdulle_heterogeneous_2012}.
The simulation of heterogeneous multiscale problems continues to be of much interest and new approaches include the work of \cite{March20, Peherstorfer20, Bastidas21}. 
Some advantages of the Equation-Free patch scheme that we adopt herein, described in \cref{sec:pat}, is their robust and efficient performance, in multiple spatial dimensions, to controllable high-order accuracy, and in the presence of micro-scale heterogeneity. 

\cite{Maclean2020b} extended the Equation-Free patch scheme (outlined in \cref{sec:dPat}) to be suitable for the efficient and accurate resolution of \emph{stationary} shocks of known location.
Here we innovatively extend this scheme in the following two critical aspects.
\begin{itemize}
    \item \cref{sec:mov,sec:movDpat} develop new techniques that automatically detect and track \emph{moving} shocks by employing moving mesh methods from, for example, \cite{Budd2009, Huang10}.
Our Equation-Free scheme needs no {a priori} information on where or when shocks form in simulation.
Instead, the scheme smoothly adapts in time the computational mesh in order to accurately simulate any number of emergent moving shocks. 
\Cref{fig:hetsin}(b) shows one example of our novel patch scheme detecting and simulating the formation of a shock, that here happens to be stationary, for a strongly heterogeneous micro-scale.
    \item \cref{sec:cons} establishes consistency results for the computational scheme for 1D~spatial systems like~\cref{hetModel}.
These results begin to establish theoretical support for both the stationary Equation-Free scheme of \cite{Maclean2020b} and the adaptive scheme introduced by \cref{sec:movPat}.
\end{itemize}
The review by \cite{dellIsola2016} challenges researchers ``to improve the capability of tailoring macroscopic homogenized models to the description of microscopic complexity''.
By accurately linking micro-scale systems with shocks to their emergent macro-scale, the developments herein provide a powerful computational capability to meet this challenge.

\Cref{sec:e1,sec:e2,sec:e3} detail three examples of our novel Equation-Free, moving and merging, patch scheme.
These examples make accurate predictions when compared to the corresponding, computationally expensive, full-domain simulation of~\cref{hetModel}.
This predictive accuracy is achieved despite computing on only a small fraction of the physical domain. 

This article focuses on simulating the formation and movement of shocks.  Developing the efficient patch scheme further to cater for the very long-time `evaporation' of dissipating shocks is left for future research.

\section{Background: \Pat{}es resolve stationary shocks} 
\label{sec:background}

This section summarises the extant patch scheme in 1D~space \cite[e.g.]{Gear03, Hyman2005, Samaey03b, Samaey04, Roberts06d}, and the scheme's extension that \cite{Maclean2020b} developed to incorporate and accurately resolve stationary shocks.
\Cref{sec:mov} then details our novel extension to moving and merging the patches in order to efficiently and accurately resolve generally emergent and moving shocks.

Generally we use uppercase symbols to denote macro-scale quantities, and lowercase to denote micro-scale quantities.

\subsection{Patch schemes}
\label{sec:pat}

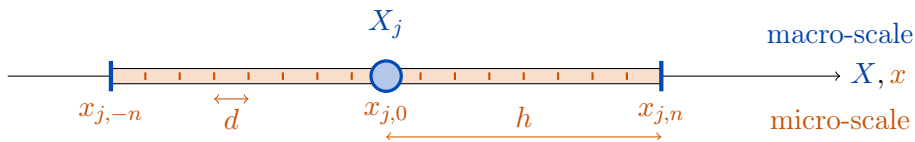
\begin{figure}
    \centering
    \tikzsetfigurename{figs/patchsetup}%
\begin{tikzpicture}
\draw[->] (0,1) -- (11,1) node[right,below=0.3cm] {\textcolor{mtxt}{micro-scale}} node[right,above=0.3cm] {\textcolor{Mtxt}{macro-scale}} node[right] {$\textcolor{Mtxt}{X},\textcolor{mtxt}{x}$};
\filldraw[fill=mcol!20,draw=black] (30/22,1.1) rectangle (10-30/22,0.9);
\foreach \x in {3,...,19} {
	\draw[mtxt,thick] (10*\x/22,1.05)--(10*\x/22,0.95);
	}
\draw[mtxt,<->] (60/22,0.7)--(70/22,0.7) node[below,midway]{$d$};
\draw[Mcol,ultra thick] (30/22,1.2)--(30/22,0.8)  node[below] {\textcolor{mtxt}{$x_{j,-n}$}}; 
\filldraw[fill=Mcol!30, draw=Mcol, very thick] (5,1) circle (0.2cm) node[above,yshift=0.35cm] {\textcolor{Mtxt}{$X_j$}} node[below,yshift=-0.2cm] {\textcolor{mtxt}{$x_{j,0}$}}; 
\draw[mcol,<->] (5,0.2)--(190/22,0.2) node[above,midway]{\textcolor{mtxt}{$h$}};
\draw[Mcol, ultra thick] (10-30/22,1.2)--(10-30/22,0.8) node[below] {\textcolor{mtxt}{$x_{j,n}$}}; 
\end{tikzpicture}
    \caption{Schematic of the micro-scale detail of a single patch (\textcolor{mtxt}{orange in the pdf version}).
Short vertical lines show the micro-grid.
Three mesh-points on the patch connect to the macro-scale (\textcolor{Mtxt}{blue in the pdf version}): the centre of the patch, and the two edge locations.
\Cref{fig:pats} shows how the field at each patch centre is interpolated between patches to provide patch-edge values that couple the patches together.}
    \label{fig:patset}
\end{figure}

\begin{figure}
    \centering
    \tikzsetfigurename{figs/patches}%
\begin{tikzpicture}
\let\scriptstyle\relax
\tikzset{->-/.style={decoration={
  markings,
  mark=at position .45 with {\arrow{>}}},postaction={decorate}}}
  
\def\h{\textwidth/70};  
\def\hh{\h * 12}; 
\def\n{3}; 
\def\xx{\h * 30}; 
\draw[->](-\xx - 6*\h/2,1)--(\xx + 6*\h/2,1)node[right]{$X$};
\foreach \x in {-2,...,2} {\ifthenelse{\NOT \x=0}{
	\filldraw[mcol!20,draw=mtxt] (\x * \hh - \h * \n - \h/2,0.92)rectangle (\x * \hh + \h * \n + \h/2,1.08);
	\draw[Mcol,ultra thick] (\x * \hh - \h * \n - \h/2,1.15)--(\x * \hh - \h * \n - \h/2,0.85);
	\draw[Mcol,ultra thick] (\x * \hh + \h * \n + \h/2,1.15)--(\x * \hh + \h * \n + \h/2,0.85);
	\draw[excol,->-] (\hh*\x,1) .. controls ($(\hh*\x/2 - \h * \n - \h/4,1.5) + abs(\x)*(2*\h,1.5*\h)$) .. (- \h * \n - \h/2,1) ;
	\filldraw[fill=Mcol!30, draw=Mcol, thick] (\hh*\x,1) circle (0.9ex) node[below=1ex]{\textcolor{Mtxt}{$\scriptstyle X_{j\ifthenelse{0>\x}{\x}{+\x}}$}}
	 node[above=0.8ex]{{\textcolor{extxt}{$\scriptstyle U_{j\ifthenelse{0>\x}{\x}{+\x}}$}}}; 
	}{
	\filldraw[mcol!20,draw=mtxt] (\x * \hh - \h * \n - \h/2,0.92)rectangle (\x * \hh + \h * \n + \h/2,1.08);
	\draw[Mcol,very thick] (\x * \hh - \h * \n - \h/2,1.15)--(\x * \hh - \h * \n - \h/2,0.85);
	\draw[Mcol,very thick] (\x * \hh + \h * \n + \h/2,1.15)--(\x * \hh + \h * \n + \h/2,0.85);
	\draw[excol,->-] (\hh*\x,1) .. controls ( \hh*\x/2 - \h * \n/2 - \h/4,1.5) .. (- \h * \n - \h/2,1) ;
	 }
	 } 
	 \def\x{0} 
	 \filldraw[fill=Mcol!30, draw=Mcol, thick] (\hh*\x,1) circle (0.9ex) node[below=1ex]{\textcolor{Mtxt}{$\scriptstyle X_{j}$}}
	 node[above=2ex]{{\textcolor{extxt}{$\scriptstyle U_j$}}}; 
\draw[Mcol,<->,thick](0, 0)--(\hh, 0) node [midway, above] {\textcolor{Mtxt}{$H_j$}};			
\end{tikzpicture}
    \caption{Schematic showing inter-patch coupling to obtain the edge value at~$x=\textcolor{Mtxt}{X_j-h}$ (\cref{fig:patset} shows the internal structure of each patch).
Interpolating centre-patch values~$\textcolor{extxt}{U_j}$ of neighbouring macro-scale nodes~$\textcolor{Mtxt}{X_j}$ determines each edge value.
This figure illustrates the case of using \(\Gamma=2\) neighbouring patches in each direction to result in a quartic interpolation.
}
    \label{fig:pats}
\end{figure}
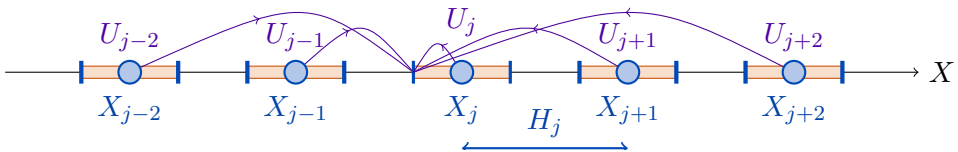

In the 1D spatial domain, set macro-scale nodes at $x=X_j$ on which the macro-scale solution is to be resolved, with spacing $H_j:=X_{j+1}- X_j$, for integer macro-scale index~\(j\) (\cref{fig:pats}).
Centred on each of these macro-scale nodes we place a spatial \emph{patch} formed by a  micro-scale mesh, a micro-grid, of $2n+1$~points with micro-scale equi-spacing~$\dx$ (\cref{fig:patset}).
The micro-scale simulation (e.g.,~\cref{hetModel}) is then only computed at the \(2n-1\)~interior points of each and every patch (\cref{fig:hetsin}(b)). 
Denote the location of the patches' micro-grid points by $x_{j,i} := X_j+\dx i$ for micro-scale index $\mathcode`\,="213B i\in\{-n,\ldots,-1,0,1,\ldots,n\}$, and the patch half-width \text{by $h:=n\dx$\,.}

\Cref{fig:patset} illustrates these details and shows three important locations for macro-scale considerations: the patch centre \emph{macro-scale node} $x_{j,0}= X_j$\,, and the two patch edge nodes at $x_{j,\pm n} = X_j \pm h$\,.
We define the patch-centre field values to be the macro-scale field values \({U_j(t)}:=u_{j,0}\) (\cref{fig:pats}) that we choose to serve as the emergent macro-scale order parameter of the multiscale patch scheme.
Many people prefer to use local averages as their macro-scale variables.
However, since the emergent macro-scale dynamics are a slow manifold of the multiscale system, and one can parameterise such a manifold in almost any way one likes \cite[\S5.3]{Roberts2014a}, then here we choose the centre-patch value as this choice most simply empowers providing \text{patch-edge values.}%
\footnote{Further, the Whitney Embedding Theorem assures us that any such choice will generally capture emergent attractors of dimension~\(N/2\) as a smooth parametrisation.}

The macro-scale values~${U_j}$ are coupled to provide field values~$u_{j,\pm n}$ at the patch edges~$x_{j,\pm n}$ (\cref{fig:pats}), closing the patch simulation.
For the edges of patch~$j$, we use polynomial interpolation of the centre-patch field values~${U_i}$ in, usually, the~$2\Gamma+1$ nearest patches to the \(j\)th~patch, for some coupling order~$\Gamma$.
\Cref{fig:pats} schematically illustrates this interpolation for one of the patch edges. 
For flexibility, denote by~\(\NN_j\) the set of neighbours of patch~\(j\) that are used in the interpolation: usually \(\NN_j:=\{j-\Gamma,\ldots,j+\Gamma\}\).
We choose to determine the edge-patch values, at each time, by the flexible classic Lagrange interpolation formula
\begin{align} \label{patCoup}
u_{j,\pm n} := \sum_{k\in\NN_j} \left(\prod_{\ell\in\NN_j\backslash \{k\}} \frac{ x_{j,\pm n} - X_\ell}{X_k - X_\ell}\right) U_k \,.
\end{align}
Other expressions for Lagrange interpolation, such as those using centred difference and mean operators \cite[]{Roberts06d, Roberts2011a, Cao2014a}, are equivalent and may be used instead.
Such Lagrange interpolation has been proven \cite[e.g.]{Roberts00a} to give a patch scheme with macro-scale predictions generally having errors~\Ord{H^{2\Gamma}} for patch spacing~\(H\).
That is, the macro-scale predictions of the patch scheme have \text{controllable errors.}

Consider the multiscale scenario where you only want to resolve macro-scale structures on length scales longer than a few times~\(H\), say, \emph{and} where the characteristic micro-scale length scales are much smaller than the macro-scale, $\dx<h\ll H$\,.
Then the patch scheme computes the micro-scale system on only the small fraction~$2h/H$ of the 1D spatial domain (\cref{fig:hetsin}(b)).
Hence, in this multiscale scenario the patch scheme is much more efficient than computing over the full spatial domain.

However, macro-scale shocks require special treatment in order to preserve accuracy and efficiency in the patch scheme.

\subsection{Shocks separate patch simulations into coupled macro-scale systems}
\label{sec:dPat}

\begin{figure}
    \centering
    \tikzsetfigurename{figs/doublepatches}%
\begin{tikzpicture}
\let\scriptstyle\relax
\tikzset{->-/.style={decoration={
  markings,
  mark=at position .45 with {\arrow{>}}},postaction={decorate}}}
  
\def\h{\textwidth/70};  
\def\hh{\h * 12}; 
\def\n{3}; 
\def\xx{\h * 30}; 
\draw[->](-\xx - 3*\h,1)--(\xx + 3*\h,1)node[right]{$x$};
\foreach \x in {-2,-1,1,2} {
	\filldraw[mcol!20,draw=mtxt] ($  (\x * \hh - \h * \n - \h/2,0.92) + \x/abs(\x)*(\n*\h,0) $)rectangle ($ (\x * \hh + \h * \n + \h/2,1.08) + \x/abs(\x)*(\n*\h,0) $);
	\draw[Mcol,ultra thick] ($  (\x * \hh - \h * \n - \h/2,1.15) + \x/abs(\x)*(\n*\h,0) $)--($  (\x * \hh - \h * \n - \h/2,0.85) + \x/abs(\x)*(\n*\h,0) $);
	\draw[Mcol,ultra thick] ($ (\x * \hh + \h * \n + \h/2,1.15) + \x/abs(\x)*(\n*\h,0) $)--($ (\x * \hh + \h * \n + \h/2,0.85) + \x/abs(\x)*(\n*\h,0) $);
	\ifthenelse{0>\x}{
	\draw[excol,->-] ($ (\hh*\x,1) + \x/abs(\x)*(\n*\h,0) $) .. controls ($(\hh*\x/2 - \h * \n/2 - \h/4,1.5) + abs(\x)*(0,1.5*\h) + \x/abs(\x)*(\n*\h,0) $) .. (- 2*\h * \n - \h/2,1) ;
	}{}
	\filldraw[fill=Mcol!30, draw=Mcol, thick] ($ (\hh*\x,1) + \x/abs(\x)*(\n*\h,0) $) circle (1ex) node[below=1.5ex]{\textcolor{Mtxt}{$\scriptstyle X_{s\ifthenelse{0>\x}{\x}{+\x}}$}}
	 node[above=1.2ex]{{\textcolor{extxt}{$\scriptstyle U_{s\ifthenelse{0>\x}{\x}{+\x}}$}}}; 
	}
\def\x{0} 
	\filldraw[fill=mcol!30, draw=mtxt] (\x * \hh - 2*\h * \n - \h/2,0.92)rectangle (\x * \hh + 2*\h * \n + \h/2,1.08);
	\draw[Mcol,very thick] (\x * \hh - 2*\h * \n - \h/2,1.15)--(\x * \hh - 2*\h * \n - \h/2,0.85);
	\draw[Mcol,very thick] (\x * \hh + 2*\h * \n + \h/2,1.15)--(\x * \hh + 2*\h * \n + \h/2,0.85);
	\draw[excol,->-] (\hh*\x-\n*\h,1) .. controls ( \hh*\x/2 - 3*\h * \n/2 - \h/4,1.5) .. (- 2*\h * \n - \h/2,1) ;
	\filldraw[fill=Mcol!30, draw=Mcol, thick] (\hh*\x-\n*\h,1) circle (1ex) node[below=1.5ex]{\textcolor{Mtxt}{$\scriptstyle X^l_{s}$}}
	 node[above=1.2ex]{{\textcolor{extxt}{$\scriptstyle U^l_s$}}}; 
	 \filldraw[fill=Mcol!30, draw=Mcol, thick] (\hh*\x+\n*\h,1) circle (1ex) node[below=1.5ex]{\textcolor{Mtxt}{$\scriptstyle X^r_{s}$}}
	 node[above=1.2ex]{{\textcolor{extxt}{$\scriptstyle U^r_s$}}}; 
\draw[Mcol,<->,thick](-\hh-\n*\h, -0.1)--(0-\n*\h, -0.1) node [midway, above] {$H_{s-1}$};			
\draw[Mcol,<->,thick](\hh+\n*\h, -0.1)--(0+\n*\h, -0.1) node [midway, above] {$H_{s}$};			
\end{tikzpicture}
\caption{A \Pat, with index \(j=s\), is drawn in the centre.
Its two macro-scale nodes are labeled \textcolor{Mtxt}{$X^l_s$} (left) and~\textcolor{Mtxt}{$X^r_s$} (right). 
Inter-patch coupling is by Lagrange interpolation~\eqref{patCoup} (\cref{fig:pats}), but the interpolation is adjusted near the \Pat, made asymmetric, so that the interpolants do not `cross' the \Pat.
For example, as indicated, we may obtain the \Pat\ left-edge value by a quadratic interpolation of~\(\textcolor{extxt}{U_{s-2},U_{s-1},U^l_s}\)}
    \label{fig:dpats}
\end{figure}
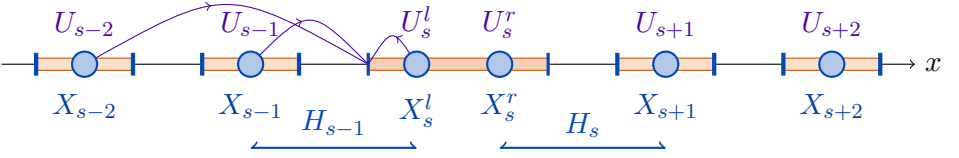

We previously captured stationary shocks in patch simulations by developing a new type of patch \cite[]{Maclean2020b}.
A so-called \emph{\Pat}, when placed over a shock in the macro-scale field, prevents otherwise ruinous interpolation between patches on either side of the shock. 
A \Pat\ is so named because its size is larger than the other patches---large enough to resolve the shock---but generally smaller than the size of other macro-scale structures the simulation aims to resolve (\cref{fig:hetsin}(b) for \(t\gtrsim1\)).
We define a \Pat\ to have two macro-scale nodes (\cref{fig:dpats}), one on either side of the micro-scale transition layer---a transition layer that appears on the macro-scale as the shock.
When the \Pat\ has patch index~\(s\), then \(X_s^r\)~denotes the location of the right node, whereas \(X_s^l\)~denotes the left (\cref{fig:dpats}).
These \Pat\ nodes are treated as spatial boundaries from the perspective of inter-patch coupling.
The overall effect is to separate a simulation into two coupled macro-scale domains, with the special \Pat\ coupling the different macro-domains by providing the information to correctly evolve both the `left' and the `right' macro-domain on either side \text{of the emergent shock.}

The macro-scale inter-patch interpolation that couples patches is treated differently in the vicinity of the \Pat\ in that the neighbour sets~\(\NN_j\) are changed.
One may use any appropriate interpolation of the~$U_j$ that treats~$U_s^l$ and~$U_s^r$ as boundary values on the left and right sides of the shock, respectively. 
With \Pat{}es there is a difference between the interpolation for the left and right edges, so for ordinary patches denote \(\NN_j^l:=\NN_j^r:=\NN_j\).
Then, for example, an interpolation with fixed bandwidth across the domain is obtained by choosing, for every edge \(e\in\{l,r\}\), 
\begin{equation*}
{\NN^e_j}:=\begin{cases}
\{j-\Gamma,\ldots,\min(j+\Gamma,s)\}
&\text{for }j<s\text{ or }(j,e)=(s,l),\\
\{\max(j-\Gamma,s),\ldots,j+\Gamma\}
&\text{for }j>s\text{ or }(j,e)=(s,r).
\end{cases}
\end{equation*}
As discussed by \cite{Maclean2020b} [\S{}IV], this choice of neighbours maintains a constant bandwidth for the interpolation across the entire spatial domain. 
Classical results \cite[e.g.]{Shoosmith75, Beyn79} assure us that despite this constant bandwidth interpolation reducing the local order of consistency with the micro-scale system, it does not affect the global \text{order of accuracy.}

Our previous research \cite[]{Maclean2020b} was successful at accurately and efficiently simulating stationary shocks of known location.
But to extend the scheme to resolve emergent and moving shocks we need to provide a mechanism to move the patches adaptively, and also to form a \Pat\ large enough to resolve the transition layers on either side of a shock. 
The next sections resolve these issues.

\section{Adaptively move patches to best resolve the macro-scale}
\label{sec:mov}

To move each patch, we implement on the macro-scale the methodology of adaptive moving meshes \cite[e.g.]{Li1998, MacKenzie07, Budd2009, Huang10}.
That is, the patch-centre locations~$X_j(t)$ are to vary in time according to a standard moving mesh method (\cref{fig:hetsin}(b)).
But within each patch we hold the micro-grid fixed---the idea is to empower a user to provide a given micro-scale code that computes the dynamics within a patch, and to wrap a flexible dynamic patch scheme around such given code.

The patch scheme and moving mesh frameworks synergise beyond the regular capability of either.
Moving patches enable patch schemes to automatically track important macro-scale features, such as strong curvatures or shocks\slash cracks. 
The patch scheme empowers us to apply the moving mesh methods to problems that are very expensive to simulate on the micro-scale: the primary example~\eqref{hetModel} considered here involves the micro-scale physics of heterogeneous diffusion and heterogeneous nonlinear advection.
Potential applications also include simulation on a complicated micro-scale domain with boundaries difficult to satisfy for moving meshes.
In addition, as confirmed by the numerical experiments of \cref{sec:e1,sec:e2,sec:e3}, the patch equations of motion may be orders of magnitude less stiff than the methods of standard moving meshes.

Before providing details on the precise formulation of moving mesh methods employed, let us motivate the discussion of \cref{sec:movDpat} on merging patches.
Applying moving mesh methods, meant for infinitesimal mesh points, to patches possessing finite widths implies patches must occasionally `collide' with each other  (\cref{fig:hetsin}(b) for \(t\approx1\)).
We exploit this feature, which typically only occurs in regions where the macro-scale field~$U_j(t)$ exhibits strong curvature.
Whenever two patches collide we assume that the macro-scale field represented by those patches exhibits some localised irregularity, such as an emergent shock.
Thus, when patches collide, we merge the patches, as detailed by \cref{sec:movDpat}, forming a \Pat\ as defined in \cref{sec:dPat}.
The `moving and merging' patch scheme embeds moving mesh methods for shock resolution in an Equation-Free framework for multiscale or multi-physics problems.
\Cref{sec:movPat} describes a separate mechanism that moves \Pat{}es.

\subsection{Moving mesh methods} \label{sec:stdMov}
On the macro-scale we need to approximate a macro-scale field~$u(x)$ with field values sampled on a mesh, not necessarily uniform.
Given a finite number of mesh points, the approach is to concentrate them on rapidly varying or complicated structures in~$u$.
One may distribute the mesh points~$X_j$ so that
\begin{align} \label{mdist}
\int_{X_j}^{X_{j+1}} \rho\big[u(x,t)\big] \d x \quad\text{is constant, independent of index }j,
\end{align}
for some chosen appropriate \emph{mesh density functional}~$\rho\big[ u(x,t) \big]$.
For example, one could better resolve regions where the solution varies rapidly in space by choosing $\rho := |\partial u/\partial x|$.
Such a distribution suffers various drawbacks when simulating time-dependent fields~$u(x,t)$: the optimal mesh may not change smoothly from one time to the next;  and the combination of a~\pde\slash\ode{}s with the mesh distribution equation~\eqref{mdist} forms a system of Differential-Algebraic Equations.
To overcome these drawbacks the moving mesh framework augments the system state~$u$ with a dynamic equation for the moving \text{mesh points~$X_j(t)$.}

Let's introduce the scheme in the context of \pde{}s on continuous space~\(x\).
The analogue of the mesh index~\(j\) is a space parameter~\(\xi\), typically chosen so that \(\xi\in[0,1]\) is diffeomorphic to the spatial domain \(x\in[a,b]\).
Then instead of describing a field as~\(u(x,t)\), we seek the field parametrically as~\(u=\hat u(\xi,t)\) at spatial location~\(x(\xi,t)\).
Then \pde{}s for the dynamics are the following.
Firstly, since $\hat u(\xi,t) = u(x(\xi,t)\C t)$, differentiation gives 
\begin{align}\label{chainU}
\left(\D t{\hat u}\right)_\xi
=\left(\D t{u}\right)_x
+\left(\D x{u}\right)_t\left(\D t{x}\right)_\xi
=\left(\D t{u}\right)_x
+\frac{\hat u_\xi}{x_\xi}\left(\D t{x}\right)_\xi
\,.
\end{align}
On the right-hand side \(\left(\D t{u}\right)_x\) is computed from the user given prescribed \pde\ (or \ode{}s such as~\eqref{hetModel}), and \(\left(\D x{u}\right)_t\) is computed as~\(\hat u_\xi/x_\xi\)\,.
Secondly, we choose the dynamics for~$x(\xi,t)$, namely \(\left(\D t{x}\right)_\xi\)\,, by some \emph{chosen moving mesh \pde}, here
\begin{align} \label{mmpde}
\D tx = \frac1{\rho[u]\, \tau} \D{\xi}{~}\left(\rho[u]\D \xi x \right) .
\end{align}
The equilibria of~\cref{mmpde} satisfy~\cref{mdist}. 
The parameter~$\tau$ in~\cref{mmpde} imposes a chosen time scale for the decay to these equilibria, and hence controls how closely the system tracks~\cref{mdist} when evolving dynamically. 

In the examples of \cref{sec:e1,sec:e2,sec:e3} we employ a macro-scale discretisation of the continuum mesh density function
\begin{equation}
\rho[u] := \left( 1 + \frac{1}{\alpha}  u_{xx}^2 \right)^{{1}/{3}} \text{,\quad where } \alpha := \left[ \frac{1}{b-a} \int_a^b \left| u_{xx}^2 \right|^{2/3} \d x\right]^3
\label{eq:rhoeg}
\end{equation}
is a regularization parameter.
\cite{Huang10} [\S2.4.4] report that this particular choice of~$\rho[u]$ is optimal when errors are measured by the $H^1$~semi-norm, and when the solution field~$u$ is obtained by linear interpolation of the field values at grid points.
The framework of moving mesh \pde s can simulate forming, moving and merging shocks---provided the micro-scale structures are smoothly resolved \cite[e.g.]{Budd2009, Huang10}.

\subsection{Move patches adaptively} 
\label{sec:movPat}
The attractively straightforward approach of this section is to move patches according to the moving mesh \pde~\eqref{mmpde} calculated using only the macro-scale nodes of each patch (\cref{fig:hetsin}(b)).
The entirety of each patch moves at one speed so that the micro-scale grid does not distort.
The patch speeds are calculated by assessing, via~\eqref{eq:rhoeg}, the macro-scale variations among neighbouring patches.

Recall that \Cref{fig:pats} shows the macro-scale view of the patch scheme.
There exists a set of macro-scale mesh points~$X_j(t)$ at the centre of each patch with associated field values~$U_j(t)$.
Suppose there are \(N\)~patches in the spatial domain~\([a,b]\).
Then for every~\(j\) indexing a patch that is interior to the domain, and recalling that \(H_j:=X_{j+1}-X_j\)\,, we discretise the moving mesh \pde~\cref{mmpde} for node locations~$X_j$ with field \text{values~$U_j$ via}
\begin{subequations} \label{mmDisc}
\allowdisplaybreaks
\begin{align} 
\label{mmX}
\de t{X_j} &= \frac{(N-1)^2}{2\rho_j \tau }\left[ (\rho_{j+1}+\rho_j) H_j - (\rho_j + \rho_{j-1}) H_{j-1} \right], \\
\label{mmR}
\rho_j &:= \left(1 + \frac{1}{\alpha} {U''_{j}}^2 \right)^{1/3},\\
\label{mmAl}
\alpha &:= \max\left\{ 1\C \left[\frac{1}{b-a}\sum_{j}  H_{j-1}
\frac12 \left({U''_{j}}^{2/3} + {U''_{j-1}}^{2/3}\right)\right]^3 \right\},
\\
U''_{j} &:= \frac2{H_j+H_{j-1}}\left[ \frac{U_{j+1} - U_j}{H_j} - \frac{U_{j} - U_{j-1}}{H_{j-1}} \right] ,
\end{align}
where the last of these approximates the second order derivative~\(U_{xx}\) at~\(X_j(t)\).
\end{subequations}
For each patch, we move all micro-grid points according to the speed of the macro-scale node mesh point: $\de t{x_{j,i} = \de t{X_j}}$ \text{for every~\(i,j\).}

If the macro-scale domain is spatially periodic, then~\eqref{mmDisc} is complete.
Alternatively, when the physical system has boundaries, for whatever micro-scale condition a user codes, then we may fix patches with index \(j=1,N\) at the boundaries to resolve any boundary layers, and use~\eqref{mmDisc} for movement of the interior patches \(j=2,\ldots,N-1\)\,.

As patches have finite widths, in a dynamically evolving simulation some patches are likely to collide in regions where the mesh density function~$\rho[u]$ is large.
We view such collision as a signal that we must resolve the emergence of some significant localised structure. 
Accordingly, we do not prevent patches from colliding, but rather introduce the procedure of the following \cref{sec:movDpat} to automatically merge a pair of colliding patches and form them into a moving {\Pat} as introduced for shocks by \cref{sec:dPat}.

\subsection{Example: accurately simulate with moving \Pat{}es}
\label{sec:e1}

This example, shown as \cref{fig:hetsin}, demonstrates that moving patches successfully gather together to resolve the formation of a shock, and then merge to form a \Pat\ on-the-fly (\cref{sec:movDpat}), without {a priori} knowledge of the location of the shock.
This example also shows that our moving patch scheme may be better than contemporary deep machine learning techniques.

The micro-scale model~\cref{hetModel} is completed by specifying the micro-grid spacing, here $\dx=0.0016$, the (\(\kappa\)-periodic) heterogeneities~$\epsilon_k$ and~$\gamma_k$, and requires initial and boundary conditions.
Herein, the heterogeneities~$\epsilon_k$ and~$\gamma_k$ are sampled from log-normal distributions,
\(\epsilon_k \sim \exp\left[\cN(0,\sige^2)\right]\) and 
\(\gamma_k \sim \exp\left[\cN(0,\sigg^2)\right]\),
and then normalised so the the mean of~$\gamma_k$ is one, and the harmonic mean of~$\epsilon_k$ is small, herein the mean is $\epsilon=0.01$.
That is, scalars~$\sige$ and~$\sigg$ characterise the magnitude of the heterogeneous variations, and~$\epsilon$ characterises a homogenised diffusion.
In this example we set the heterogeneous strengths $\sige=1$ and $\sigg=2$ with micro-grid periodicity $\kappa=5$, so that the micro-scale, with its log-normal coefficients, is strongly heterogeneous.
\cref{tblone} lists the specific realisation used in this example.
\begin{table}
\caption{\label{tblone}specific heterogeneous coefficients, five-periodic, for the \ode{}s~\eqref{hetModel} used in the example of \protect\cref{sec:e1,fig:hetsin,fig:sinhetDetails}.}
\begin{equation*}
\begin{array}{r|lllll}
\gamma_k& 0.38  &  1.36   & 0.63 &   3.97   & 0.19
\\
\epsilon_k&  0.003  &  0.033  &  0.14   & 0.018   & 0.012 
\end{array}
\end{equation*}
\end{table}

We choose the spatial domain to be~\([-\pi,\pi]\), initial conditions $u(x,0) = -\sin(x)$ sampled on the micro-grid, so that the full domain simulation of \cref{fig:hetsin}(a) has \(4\,000\)~micro-grid points in space, and choose boundary conditions $u(\pm\pi,t)=0$.
The computationally expensive full domain simulation \cref{fig:hetsin}(a) is accurately predicted by our more efficient moving patch scheme, \cref{fig:hetsin}(b).

The moving patch scheme requires that we specify the patch configuration, and parametrise how they move.
In this example we initially distribute $N=26$~patches, evenly spaced, each consisting of $50$~micro-scale points, so \(n=25\), as shown in \cref{fig:hetsin}(b).
Patches are coupled together via interpolation of order $\Gamma=6$ (\cref{sec:pat}).
Here, initially there is no patch covering $x=0$, where the shock emerges, but that this is no obstacle to \text{accurate simulation.}

The parameter~\(\tau\) controls patch movement:  it specifies the time-scale on which patches approach a quasi-equilibria distribution of~\cref{mdist}.
Setting this time-scale very small means that patches are very nearly perfectly placed at all times. 
However, in practice small~\(\tau\) induces numerical error in the simulation, and markedly increases the stiffness of the dynamical system.
For standard adaptive meshes, \cite{Huang10} suggest that one should choose the stiffness parameter $\tau\approx0.01$.
However, we discovered that adaptively moving patches can be much less stiff than moving meshes, and consequently choices $\tau\in\{0.1,1,10\}$ are viable (that is, the adaptive moving patch scheme may be a hundred, or possibly a thousand, times less stiff than an adaptive mesh scheme).
The mesh redistribution parameter for patches is~$\tau=10$ (\cref{sec:stdMov}), a thousand times less stiff than recommended for standard moving meshes \cite[]{Huang10}.

\Cref{fig:hetsin}(b) shows how the moving patches congregate to resolve the shock that forms in this example at about \((x,t)=(0,1)\).
However, the shock itself, albeit small on the macro-scale (\cref{fig:hetsin}(a)), is a relatively large structure on the micro-scale. 
\Cref{fig:hetsin}(b) shows how six of the moving patches `collide' together at the shock. 
So, in order to resolve the shock properly on its meso-scale we need to merge patches together to form a meso-patch as introduced by \cref{sec:dPat}. 
The next \cref{sec:movDpat} details how we smoothly merge patches.
\begin{figure}
\includegraphics{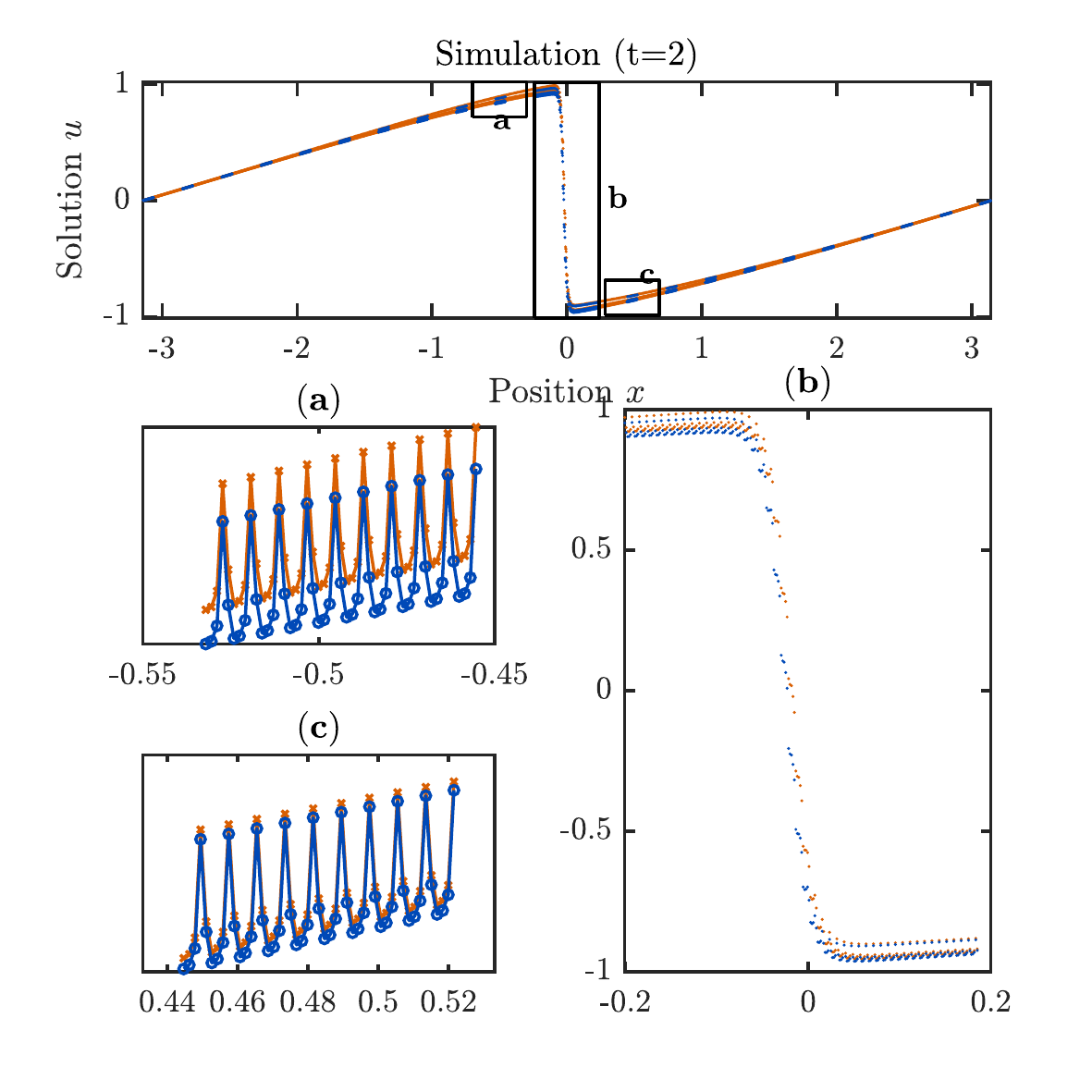}
\caption{Detailed view at the final time of the simulations shown in \cref{fig:hetsin}, described in \cref{sec:e1}.
The macro-scale view of the top graph shows the patches (blue dots) overlay the exact simulation (orange dots) so closely that they cannot be distinguished.
The insets~(a)\&(c) show micro-scale details of the two patches on either side of the \Pat: the gradient, period, amplitude, and value of the micro-scale heterogeneities are predicted almost perfectly.
Inset~(b) details the \Pat\ covering the shock.
The \Pat\ (blue dots) align with the exact solution (orange dots) on either side of the shock, with a small micro-scale phase shift of the shock in the patch prediction.}
\label{fig:sinhetDetails}
\end{figure}%
With patches merging to form a meso-patch, \Cref{fig:sinhetDetails} shows details at the final time $t=2$ of the patch predictions for this example.
The phase and periodicity of the micro-scale variations are accurately represented.
Even the location of the shock is predicted accurately with only a \text{small micro-scale shift.} 

We also confirmed (but do not include the results here) that the macroscale predictions and simulation skill of the patch scheme were similar for all phase shifts of the heterogeneous coefficients in \cref{tblone}.

We emphasise that these, and subsequent, successful predictions by the patch scheme are arrived at with \emph{no} knowledge of a macro-scale homogenisation of the nonlinear heterogeneous micro-scale system~\eqref{hetModel}, and \emph{no} knowledge of any Rankine--Hugoniot like condition for an emergent macro-scale shock. 

\paragraph{Remark on Machine Learning}
Our simulations for this example are directly comparable to some of those invoking Deep Learning.
\cite{Lu21} [\S4.2] use the homogeneous equivalent of this experiment to test a Physics-Informed Neural Network (\textsc{pinn}).
An important difference in the computational experiments is that our problem possesses strongly heterogeneous, micro-scale, advection and diffusion effects.
Otherwise the problems are the same up to scaling terms.
\cite{Lu21} also resolve shocks via a spatial grid of 2\,500~points, which then requires adaptive refinement by adding even more points, whereas our patch code uses just 1000~micro-scale points throughout.
\cite{Lu21} find that their \textsc{pinn} typically simulates with $L^2$~relative errors of about $10$--$15\%$.
We computed the maximum $L^2$~relative error over all simulated times for our patch predictions at the macro-scale nodes, and discovered that our adaptive patch scheme has $L^2$~relative errors of~$1.3\%$, about a factor of ten better than the \textsc{pinn}---and does so for a much more difficult problem of heterogeneity in both the \text{diffusion and advection.}

\section[Form \Pat{}es and selfishly move them]{Form \Pat{}es on-the-fly and `selfishly' move them to track macro-scale shocks}
\label{sec:movDpat}

This section develops mechanisms to form moving \Pat{}es on-the-fly. 
That is, firstly, shock locations no longer need to be known a priori, and secondly, the \Pat{}es now track moving shocks---shocks no longer need be stationary.
These mechanisms exploit the patch collisions that occur when patches move according to the moving mesh \pde\ described in \cref{sec:movPat} and seen in the example of \cref{sec:e1}.

With merging patches, the number of micro-grid points in each patch (previously~$2n+1$ for every patch) may differ between patches.
Indeed it may be desirable to have wider patches at domain boundaries in order to adequately resolve micro-scale boundary layers in space.
Consequently, we adjust the notation so that there are $2n_j+1$ micro-grid points in the $j$th~patch.
Correspondingly, the patch half-width $h_j:=n_j\dx$ also now depends upon the patch index~$j$.
We omit denoting the time dependence of these quantities since they are piecewise constant in time.

\subsection{A \Pat\ emerges from patch collisions} 
\label{sec:patMrg}

Suppose two patches collide at time~\(t'\) (\cref{fig:hetsin}(b) for \(t'\approx1\))---say the \(s\)th~patch and the \((s+1)\)th~patch collide.
Collision means that their edges touch with $x_{s,n_s}(t') = x_{s+1,-n_{s+1}}(t')$.
At that time in the simulation, we interrupt the simulation and update the details of the patches as follows.
The outcome is to merge the two colliding patches as a single \Pat, the new $s$th~patch, and decrement the number of patches, $N\mapsto N-1$\,.

In detail, save the collision location $x\mrg := x_{s,n_s}= x_{s+1,-n_{s+1}}$\,, 
and the average field value at the colliding patch edges, $u\mrg := \tfrac12(u_{s,n_s} + u_{s+1,-n_{s+1}})$.
The new merged \Pat\ has a larger half-width, containing micro-grid points indexed by \(i'=-n',\ldots,n'\)\, where $n':= n_{s}+n_{s+1}$. These new grid points include all the micro-grid points of the former $s$th and $(s+1)$th~patches.
The patch indices are \text{related by}
\begin{equation*}
i'=\begin{cases}
i+n_s-n'&i\in[-n_s,n_s]\text{ of \(s\)th patch}, \\
i-n_{s+1}+n'&i\in[-n_{s+1},n_{s+1}]\text{ of \((s+1)\)th patch}.
\end{cases}
\end{equation*}
This merger, at new \Pat\ index~$i'=n_s-n_{s+1}$\,, generates one additional \emph{interior} micro-grid point for which the micro-scale \ode~\eqref{hetModel} is thereafter applied---for times \(t>t'\).
Its `initial' location at time~\(t'\) is~\(x\mrg\) with `initial' field value~$u\mrg$.
The \Pat\ field values at the other $2n'-2$ micro-grid points internal to the \Pat\ are not altered from those of the two colliding patches. 
The new \Pat\ is then indexed as the new \(s\)th~patch, and the patches with higher indices have their index decremented by one. 
The result of this patch merge is always a \Pat, and in all cases the two macro-scale nodes of the new \Pat\ are uniquely specified as \text{shown in \cref{fig:mrg}.}

\begin{figure} 
\centering  
\begin{enumerate}[label=(\alph*) , ref=\cref{fig:mrg}(\alph*)]
\item \label{fig:mrgSS} Merging two ordinary patches: \\
    \tikzsetfigurename{figs/ssmerge}%
\begin{tikzpicture}
\let\scriptstyle\relax
\tikzset{->-/.style={decoration={
  markings,
  mark=at position .45 with {\arrow{>}}},postaction={decorate}}}
  
\def\h{\textwidth/52} 
\def\bl{\h} 

\def\hgt{1.5} 
\def\lhgt{\hgt+0.3} 
\def\rhgt{\hgt-0.3}
\draw[->] (0,\hgt) -- (11,\hgt)  node[right] {$\textcolor{Mtxt}{X}$} 
 node[yshift=2ex, left=\h] {time \(t'-\)}; 

\def\nl{10} 
\def\nr{10} 
\def\br{\bl+\nl*\h} 
\def\b{\bl} \def\n{\nl} 
\filldraw[fill=mcol!20,draw=black] (\b, \lhgt+0.1) rectangle (\b+\n*\h, \lhgt-0.1);
\foreach \x in {1,...,\n} {
	\draw[mtxt,thick] (\b+\x*\h,\lhgt+0.05)--(\b+\x*\h,\lhgt-0.05);
	}
\draw[mtxt,|-|] (\b+\h,\lhgt+0.3)--(\b+2*\h,\lhgt+0.3) node[above,midway]{$d$};
\draw[Mcol,ultra thick] (\b,\lhgt+0.2)--(\b,\lhgt-0.2)  ; 
\filldraw[fill=Mcol!30, draw=Mcol, very thick] (\b+5*\h,\lhgt) circle (0.13cm) node[above,yshift=+0.1cm] {\textcolor{Mtxt}{$X_s$}};
\draw[Mcol, ultra thick] (\b+\n*\h,\lhgt+0.2)--(\b+\n*\h,\lhgt-0.2) ; 

\def\b{\br}  \def\n{\nr}
\filldraw[fill=mcol!20,draw=black] (\b, \rhgt+0.1) rectangle (\b+\n*\h, \rhgt-0.1);
\foreach \x in {1,...,\n} {
	\draw[mtxt,thick] (\b+\x*\h,\rhgt+0.05)--(\b+\x*\h,\rhgt-0.05);
	}
\draw[Mcol,ultra thick] (\b,\rhgt+0.2)--(\b,\rhgt-0.2)  ; 
\filldraw[fill=Mcol!30, draw=Mcol, very thick] (\b+5*\h,\rhgt) circle (0.13cm) node[below,yshift=-0.2cm] {\textcolor{Mtxt}{$X_{s+1}$}};
\draw[Mcol, ultra thick] (\b+\n*\h,\rhgt+0.2)--(\b+\n*\h,\rhgt-0.2) ;  

\def\hgt{0} 
\def\b{\bl} 
\def\n{20} 
\draw[->] (0,\hgt) -- (11,\hgt)  node[right] {$\textcolor{Mtxt}{X}$}
 node[yshift=2ex, left=\h] {time \(t'+\)}; 
\filldraw[fill=mcol!20,draw=black] (\b, \hgt+0.1) rectangle (\b+\n*\h, \hgt-0.1);
\foreach \x in {1,...,\n} {
	\draw[mtxt,thick] (\b+\x*\h,\hgt+0.05)--(\b+\x*\h,\hgt-0.05);
	}
\draw[Mcol,ultra thick] (\b,\hgt+0.2)--(\b,\hgt-0.2)  ; 
\draw[Mcol, ultra thick] (\b+\n*\h,\hgt+0.2)--(\b+\n*\h,\hgt-0.2) ;  
\filldraw[fill=Mcol!30, draw=Mcol, very thick] (\bl+5*\h,\hgt) circle (0.13cm) node[below,yshift=-0.2cm] {\textcolor{Mtxt}{$X_s^l$}};
\filldraw[fill=Mcol!30, draw=Mcol, very thick] (\br+5*\h,\hgt) circle (0.13cm) node[below,yshift=-0.2cm] {\textcolor{Mtxt}{$X_{s}^r$}};

\end{tikzpicture}
\item \label{fig:mrgSD} Merging an ordinary patch with a \Pat:\\
    \tikzsetfigurename{figs/sdmerge}%
\begin{tikzpicture}
\let\scriptstyle\relax
\tikzset{->-/.style={decoration={
  markings,
  mark=at position .45 with {\arrow{>}}},postaction={decorate}}}
  
\def\h{\textwidth/52} 
\def\bl{\h} 

\def\hgt{1.5} 
\def\lhgt{\hgt+0.3} 
\def\rhgt{\hgt-0.3}
\draw[->] (0,\hgt) -- (11,\hgt)  node[right] {$\textcolor{Mtxt}{X}$}
 node[yshift=2ex, left=\h] {time \(t'-\)}; 

\def\nl{10} 
\def\nr{21} 
\def\br{\bl+\nl*\h} 
\def\b{\bl} \def\n{\nl} 
\filldraw[fill=mcol!20,draw=black] (\b, \lhgt+0.1) rectangle (\b+\n*\h, \lhgt-0.1);
\foreach \x in {1,...,\n} {
	\draw[mtxt,thick] (\b+\x*\h,\lhgt+0.05)--(\b+\x*\h,\lhgt-0.05);
	}
\draw[mtxt,|-|] (\b+\h,\lhgt+0.3)--(\b+2*\h,\lhgt+0.3) node[above,midway]{$d$};
\draw[Mcol,ultra thick] (\b,\lhgt+0.2)--(\b,\lhgt-0.2)  ; 
\filldraw[fill=Mcol!30, draw=Mcol, very thick] (\b+5*\h,\lhgt) circle (0.13cm) node[above,yshift=+0.1cm] {\textcolor{Mtxt}{$X_s$}};
\draw[Mcol, ultra thick] (\b+\n*\h,\lhgt+0.2)--(\b+\n*\h,\lhgt-0.2) ; 

\def\b{\br}  \def\n{\nr}
\filldraw[fill=mcol!20,draw=black] (\b, \rhgt+0.1) rectangle (\b+\n*\h, \rhgt-0.1);
\foreach \x in {1,...,\n} {
	\draw[mtxt,thick] (\b+\x*\h,\rhgt+0.05)--(\b+\x*\h,\rhgt-0.05);
	}
\draw[Mcol,ultra thick] (\b,\rhgt+0.2)--(\b,\rhgt-0.2)  ; 
\filldraw[fill=Mcol!30, draw=Mcol, very thick] (\b+5*\h,\rhgt) circle (0.13cm) node[below,yshift=-0.2cm] {\textcolor{Mtxt}{$X_{s+1}^l$}};
\filldraw[fill=Mcol!30, draw=Mcol, very thick] (\b+16*\h,\rhgt) circle (0.13cm) node[below,yshift=-0.2cm] {\textcolor{Mtxt}{$X_{s+1}^r$}};
\draw[Mcol, ultra thick] (\b+\n*\h,\rhgt+0.2)--(\b+\n*\h,\rhgt-0.2) ;  

\def\hgt{0} 
\def\b{\bl} 
\def\n{31} 
\draw[->] (0,\hgt) -- (11,\hgt)  node[right] {$\textcolor{Mtxt}{X}$}
 node[yshift=2ex, left=\h] {time \(t'+\)}; 

\filldraw[fill=mcol!20,draw=black] (\b, \hgt+0.1) rectangle (\b+\n*\h, \hgt-0.1);
\foreach \x in {1,...,\n} {
	\draw[mtxt,thick] (\b+\x*\h,\hgt+0.05)--(\b+\x*\h,\hgt-0.05);
	}
\draw[Mcol,ultra thick] (\b,\hgt+0.2)--(\b,\hgt-0.2)  ; 
\draw[Mcol, ultra thick] (\b+\n*\h,\hgt+0.2)--(\b+\n*\h,\hgt-0.2) ;  
\filldraw[fill=Mcol!30, draw=Mcol, very thick] (\bl+5*\h,\hgt) circle (0.13cm) node[below,yshift=-0.2cm] {\textcolor{Mtxt}{$X_s^l$}};
\filldraw[fill=Mcol!30, draw=Mcol, very thick] (\br+16*\h,\hgt) circle (0.13cm) node[below,yshift=-0.2cm] {\textcolor{Mtxt}{$X_{s}^r$}};

\end{tikzpicture} 
\item Merging two \Pat{}es:\\
    \tikzsetfigurename{figs/ddmerge}%
\begin{tikzpicture}
\let\scriptstyle\relax
\tikzset{->-/.style={decoration={
  markings,
  mark=at position .45 with {\arrow{>}}},postaction={decorate}}}
  
\def\h{\textwidth/52} 
\def\bl{\h} 

\def\hgt{1.5} 
\def\lhgt{\hgt+0.3} 
\def\rhgt{\hgt-0.3}
\draw[->] (0,\hgt) -- (11,\hgt)  node[right] {$\textcolor{Mtxt}{X}$}
 node[yshift=2ex, left=\h] {time \(t'-\)}; 

\def\nl{19} 
\def\nr{19} 
\def\br{\bl+\nl*\h} 
\def\b{\bl} \def\n{\nl} 
\filldraw[fill=mcol!20,draw=black] (\b, \lhgt+0.1) rectangle (\b+\n*\h, \lhgt-0.1);
\foreach \x in {1,...,\n} {
	\draw[mtxt,thick] (\b+\x*\h,\lhgt+0.05)--(\b+\x*\h,\lhgt-0.05);
	}
\draw[mtxt,|-|] (\b+\h,\lhgt+0.3)--(\b+2*\h,\lhgt+0.3) node[above,midway]{$d$};
\draw[Mcol,ultra thick] (\b,\lhgt+0.2)--(\b,\lhgt-0.2)  ; 
\filldraw[fill=Mcol!30, draw=Mcol, very thick] (\b+5*\h,\lhgt) circle (0.13cm) node[above,yshift=+0.1cm] {\textcolor{Mtxt}{$X_s^l$}};
\filldraw[fill=Mcol!30, draw=Mcol, very thick] (\b+14*\h,\lhgt) circle (0.13cm) node[above,yshift=+0.1cm] {\textcolor{Mtxt}{$X_{s}^r$}};
\draw[Mcol, ultra thick] (\b+\n*\h,\lhgt+0.2)--(\b+\n*\h,\lhgt-0.2) ; 

\def\b{\br}  \def\n{\nr}
\filldraw[fill=mcol!20,draw=black] (\b, \rhgt+0.1) rectangle (\b+\n*\h, \rhgt-0.1);
\foreach \x in {1,...,\n} {
	\draw[mtxt,thick] (\b+\x*\h,\rhgt+0.05)--(\b+\x*\h,\rhgt-0.05);
	}
\draw[Mcol,ultra thick] (\b,\rhgt+0.2)--(\b,\rhgt-0.2)  ; 
\filldraw[fill=Mcol!30, draw=Mcol, very thick] (\b+5*\h,\rhgt) circle (0.13cm) node[below,yshift=-0.2cm] {\textcolor{Mtxt}{$X_{s+1}^l$}};
\filldraw[fill=Mcol!30, draw=Mcol, very thick] (\b+14*\h,\rhgt) circle (0.13cm) node[below,yshift=-0.2cm] {\textcolor{Mtxt}{$X_{s+1}^r$}};
\draw[Mcol, ultra thick] (\b+\n*\h,\rhgt+0.2)--(\b+\n*\h,\rhgt-0.2) ;  

\def\hgt{0} 
\def\lhgt{\hgt+0.1}
\def\rhgt{\hgt-0.1}
\def\b{\bl} 
\def\n{38} 
\draw[->] (0,\hgt) -- (11,\hgt)  node[right] {$\textcolor{Mtxt}{X}$}
 node[yshift=3ex, left=\h] {time \(t'+\)}; 

\filldraw[fill=mcol!20,draw=black] (\b, \hgt+0.1) rectangle (\b+\h*\n, \hgt-0.1);
\foreach \x in {1,...,\n} {
	\draw[mtxt,thick] (\b+\x*\h,\hgt+0.05)--(\b+\x*\h,\hgt-0.05);
	}
\draw[Mcol,ultra thick] (\b,\hgt+0.2)--(\b,\hgt-0.2)  ; 
\draw[Mcol, ultra thick] (\b+\n*\h,\hgt+0.2)--(\b+\n*\h,\hgt-0.2) ;  
\filldraw[fill=Mcol!30, draw=Mcol, very thick] (\bl+5*\h,\hgt) circle (0.13cm) node[below,yshift=-0.2cm] {\textcolor{Mtxt}{$X_s^l$}};
\filldraw[fill=Mcol!30, draw=Mcol, very thick] (\br+14*\h,\hgt) circle (0.13cm) node[below,yshift=-0.2cm] {\textcolor{Mtxt}{$X_{s}^r$}};

\end{tikzpicture}
\end{enumerate}
\caption{A \Pat\ is uniquely constructed from the merger at time~\(t'\) of two patches of any types.
Each of~(a), (b), and~(c) shows two colliding patches at time~\(t'-\), and, at time~\(t'+\), the \Pat\ resulting from their merge. }
\label{fig:mrg} 
\end{figure}

Recall we interpret the collision of two patches as meaning that there is some emergent micro-scale structure, non-smooth on the macro-scale, that needs to be resolved inside the \Pat.
Consequently, the merged \Pat\ is to contain two macro-scale nodes as defined by \cref{sec:dPat}.
For conciseness, define that for an ordinary patch with one macro-scale node the locations \(X_j^l:=X_j^r:=X_j\).  
Then, as illustrated in \cref{fig:mrg}, the new \(X_s^l(t'+):=X_s^l(t'-)\) and the new \(X_s^r(t'+):=X_{s+1}^r(t'-)\).
Further, assuming some sort of shock in the middle of a \Pat, the macro-scale interpolation should not go across the \Pat, and so we adjust the neighbourhoods~\(\NN^e_j\) of interpolation (\cref{sec:pat,sec:dPat}).
Then, for every edge \(e\in\{l,r\}\), \text{after the merger}
\begin{equation*}
{\NN^e_j}':=\begin{cases}
\NN^e_j\backslash\{s+1,s+2,\ldots\}
&\text{for }j<s\text{ or }(j,e)=(s,l),\\
(\NN^e_{j+1}-1)\backslash\{s-1,s-2,\ldots\}
&\text{for }j>s\text{ or }(j,e)=(s,r).
\end{cases}
\end{equation*}
As in \cref{sec:dPat}, these choices for the neighbour sets of the interpolation preserves the bandwidth of the interpolation without affecting the global order of accuracy \cite[e.g.]{Shoosmith75, Beyn79}.

\paragraph{Preserve the cycle of heterogeneity}
In a system with micro-scale periodic heterogeneity, \cite{Bunder2013b} established that highly accurate patch simulation resulted from choosing \(n_j\) to be an integral multiple of the heterogeneity period~\(\kappa\).
That is, it is best to choose and preserve that $n_j \bmod \kappa=0$ for every patch and \Pat.

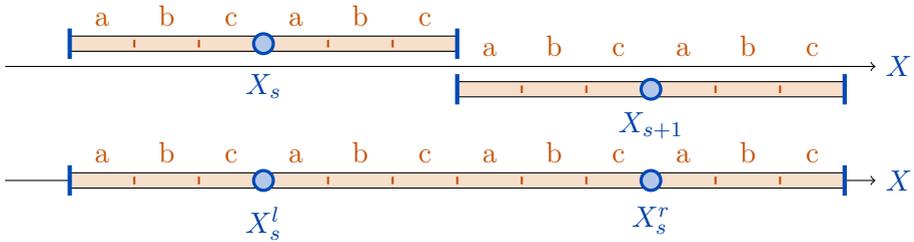
\begin{figure}\centering
    \tikzsetfigurename{figs/hetmerge}%
\begin{tikzpicture}
\let\scriptstyle\relax
\tikzset{->-/.style={decoration={
  markings,
  mark=at position .45 with {\arrow{>}}},postaction={decorate}}}
  
\def\h{\textwidth/15} 
\def\bl{\h} 

\def\hgt{1.5} 
\def\lhgt{\hgt+0.3} 
\def\rhgt{\hgt-0.3}
\draw[->] (0,\hgt) -- (11.5,\hgt)  node[right] {$\textcolor{Mtxt}{X}$};

\def\nl{6} 
\def\nr{6} 
\def\nNod{3} 
\def\br{\bl+\nl*\h} 
\def\b{\bl} \def\n{\nl} 
\filldraw[fill=mcol!20,draw=black] (\b, \lhgt+0.1) rectangle (\b+\n*\h, \lhgt-0.1);
\foreach \x in {1,...,\n} { 
	\draw[mtxt,thick] (\b+\x*\h,\lhgt+0.05)--(\b+\x*\h,\lhgt-0.05);
	}	
\foreach \x in {1,4,...,\n} { 
	\node[mtxt,above] at (\b+\x*\h-0.5*\h,\lhgt+0.1){a};
	}
\foreach \x in {2,5,...,\n} {
	\node[mtxt,above] at (\b+\x*\h-0.5*\h,\lhgt+0.1){b};
		}
\foreach \x in {3,6,...,\n} {
	\node[mtxt,above] at (\b+\x*\h-0.5*\h,\lhgt+0.1){c};
	}
%

\draw[Mcol,ultra thick] (\b,\lhgt+0.2)--(\b,\lhgt-0.2)  ; 
\filldraw[fill=Mcol!30, draw=Mcol, very thick] (\b+\nNod*\h,\lhgt) circle (0.13cm) node[below,yshift=-1.5ex] {\textcolor{Mtxt}{$X_s$}};
\draw[Mcol, ultra thick] (\b+\n*\h,\lhgt+0.2)--(\b+\n*\h,\lhgt-0.2) ; 

\def\b{\br}  \def\n{\nr}
\filldraw[fill=mcol!20,draw=black] (\b, \rhgt+0.1) rectangle (\b+\n*\h, \rhgt-0.1);
\foreach \x in {1,...,\n} {
	\draw[mtxt,thick] (\b+\x*\h,\rhgt+0.05)--(\b+\x*\h,\rhgt-0.05);
	}
\foreach \x in {1,4,...,\n} { 
	\node[mtxt,above] at (\b+\x*\h-0.5*\h,\rhgt+0.3){a};
	}
\foreach \x in {2,5,...,\n} {
	\node[mtxt,above] at (\b+\x*\h-0.5*\h,\rhgt+0.3){b};
		}
\foreach \x in {3,6,...,\n} {
	\node[mtxt,above] at (\b+\x*\h-0.5*\h,\rhgt+0.3){c};
	}
\draw[Mcol,ultra thick] (\b,\rhgt+0.2)--(\b,\rhgt-0.2)  ; 
\filldraw[fill=Mcol!30, draw=Mcol, very thick] (\b+\nNod*\h,\rhgt) circle (0.13cm) node[below,yshift=-0.15cm] {\textcolor{Mtxt}{$X_{s+1}$}};
\draw[Mcol, ultra thick] (\b+\n*\h,\rhgt+0.2)--(\b+\n*\h,\rhgt-0.2) ;  

\def\hgt{0} 
\def\b{\bl} 
\def\n{12} 
\draw[->] (0,\hgt) -- (11.5,\hgt)  node[right] {$\textcolor{Mtxt}{X}$};
\filldraw[fill=mcol!20,draw=black] (\b, \hgt+0.1) rectangle (\b+\n*\h, \hgt-0.1);
\foreach \x in {1,...,\n} {
	\draw[mtxt,thick] (\b+\x*\h,\hgt+0.05)--(\b+\x*\h,\hgt-0.05);
	}
\foreach \x in {1,4,...,\n} { 
	\node[mtxt,above] at (\b+\x*\h-0.5*\h,\lhgt-0.2){a};
	}
\foreach \x in {2,5,...,\n} {
	\node[mtxt,above] at (\b+\x*\h-0.5*\h,\lhgt-0.2){b};
		}
\foreach \x in {3,6,...,\n} {
	\node[mtxt,above] at (\b+\x*\h-0.5*\h,\lhgt-0.2){c};
	}
\draw[Mcol,ultra thick] (\b,\hgt+0.2)--(\b,\hgt-0.2)  ; 
\draw[Mcol, ultra thick] (\b+\n*\h,\hgt+0.2)--(\b+\n*\h,\hgt-0.2) ;  
\filldraw[fill=Mcol!30, draw=Mcol, very thick] (\bl+\nNod*\h,\hgt) circle (0.13cm) node[below,yshift=-0.2cm] {\textcolor{Mtxt}{$X_s^l$}};
\filldraw[fill=Mcol!30, draw=Mcol, very thick] (\br+\nNod*\h,\hgt) circle (0.13cm) node[below,yshift=-0.2cm] {\textcolor{Mtxt}{$X_{s}^r$}};

\end{tikzpicture}
\caption{Schematic for merging two patches with heterogeneous micro-scale structure of period \(\kappa=3\)\,.
Here the micro-grid labels~a, b, and~c represent the three different heterogeneous values in each period.
Crucially, merging two patches when their edges touch preserves the correct a-b-c cycle in the resultant \Pat.
This preservation of the correct micro-scale cycle of heterogeneity generalises to any period~$\kappa$ provided that it always \text{satisfies $n_j \bmod \kappa=0$\,.}}
\label{fig:hetmerge}
\end{figure}

Consequently, when two patches merge to form a \Pat, and similarly for the other merges of \cref{fig:mrg}, we need to preserve both that $n_j \bmod \kappa=0$ and also that the micro-scale periodic cycle is preserved.
As illustrated by \cref{fig:hetmerge}, we preserve both these properties by merging two patches when their two edges first touch, as discussed above.
That is, our above scheme to merge patches preserves the micro-scale periodic cycle, and also preserves that $n_j \bmod \kappa=0$ for potentially highly accurate predictions.

\subsection{Track moving macro-scale features with \Pat{}es}
\label{sectmmfwmp}
We require the \Pat{}es to cover shocks, structures which are non-smooth on the macro-scale.
But generally shocks move, thus each \Pat\ must move with the shock structure it contains.
This requirement replaces the moving mesh rules of \cref{sec:movPat} for ordinary patch movement.
Therefore the \Pat{}es must move without reference to other nearby patches (a choice referred to as `selfish movement') since resolving its shock is the \text{role of a \Pat.}

We aim to move the \Pat\ so that the shock is reasonably centrally located in the patch.
This movement of a \Pat\ replaces all Rankine--Hugoniot like conditions because the moving patch scheme is to apply to complex systems where Rankine--Hugoniot conditions are unknown.
In essence, the left\slash right side information used in Rankine--Hugoniot conditions are replaced by the \emph{separate} interpolations used on the left\slash right edges of a \Pat.
Possible future research could show, for some basic systems, that our \Pat\ rules would effectively reproduce the appropriate Rankine--Hugoniot  conditions.

Here we need each \Pat\ to interrogate the micro-scale simulation on its internal micro-grid to discern the approximate location of a shock, or other complicated micro-scale structure that emerges on the macro-scale.
Suppose the $j$th~patch is a \Pat.
The centre location of the \Pat\ is~\(x_{j,0}\). 
We conjecture that the most important micro-scale structures are where micro-scale gradients are steep, that is, $(\D xu)^2$ is large.
We estimate the location of steep gradients in the \Pat\ by choosing a weighting of positions by~$(\D xu)^2$.
Hence we define the location of interesting steep gradients in the \Pat\ as
\begin{subequations}\label{eqs:microDblMov}%
\begin{equation}
\hat X_j(t) :=   \frac{\int_{x_{j,-n}}^{x_{j,n}}  x \left(\D xu\right)^2\d x }{\int_{x_{j,-n}}^{x_{j,n}} \left(\D xu\right)^2\d x } 
\approx   \frac{\sum_{i=-n}^{n-\kappa}  \tfrac{x_{i+\kappa}+x_{i}}{2} \left(\frac{u_{i+\kappa}-u_{i}}{x_{i+\kappa}-x_{i}}\right)^2}{\sum_{i=-n}^{n-\kappa}  \left(\frac{u_{i+\kappa}-u_{i}}{x_{i+\kappa}-x_{i}}\right)^2}
\label{eq:hatX}
\end{equation}
for micro-scale periodicity~\(\kappa\) (\(\kappa=1\) if homogeneous), and where all variables in the right-hand side are implicitly for the \(j\)th~patch.
Then we move the \Pat\ as a whole to bring its central location~\(x_{j,0}(t)\) towards the important micro-scale structure, as estimated by~\(\hat X_j(t)\), according to
\begin{equation}\label{microDblMov}
\frac{dx_{j,i}}{dt} = \frac1\beta \left( \hat X_j - x_{j,0} \right),
\quad \text{for every }i,
\end{equation}
\end{subequations}
where $\beta$ is a chosen time-scale.
The parameter~$\beta$ specifies the time-scale on which \Pat{}es move towards the quasi-equilibria of~\eqref{microDblMov} in order to cover the emergent shocks, or other such localised meso-scale structures.
Our computational experiments indicate that parameter~$\beta$ should be chosen so that~$1/\beta$ is larger than the rate at which shocks move.
For the examples of \cref{sec:e2,sec:e3} we \text{choose $\beta=1$\,.}

When we simulate shocks with micro-scale heterogeneities (e.g., \cref{sec:e1,sec:e2,sec:e3}),  the micro-scale structures (e.g., insets of \cref{fig:hetsin}) obscure meso-scale shocks by the induced heterogeneous gradients.
This particularly occurs when there are heterogeneous advection terms in the micro-scale system.
For such problems we smooth the summands in the numerator and denominator of \cref{eq:hatX} by using difference estimators over one micro-period~\(\kappa\) to estimate a meso-scale~\(\D xu\). 
Alternatively, one could use a rolling mean~\(\bar u_i\) instead of the \text{micro-grid value~\(u_i\).}

There appear to be many alternatives to our choice of~\cref{eqs:microDblMov}: the best option for a particular system may depend on physical details of the micro-scale modelling.
Nonetheless, the above gives a general approach that should be useful for many \text{micro-scale systems.}

\subsection{Example: moving and merging \Pat{}es}
\label{sec:e2}

\begin{figure}
\centering
\begin{tabular}{c c}
(a)&(b)\\
\includegraphics[scale=0.95]{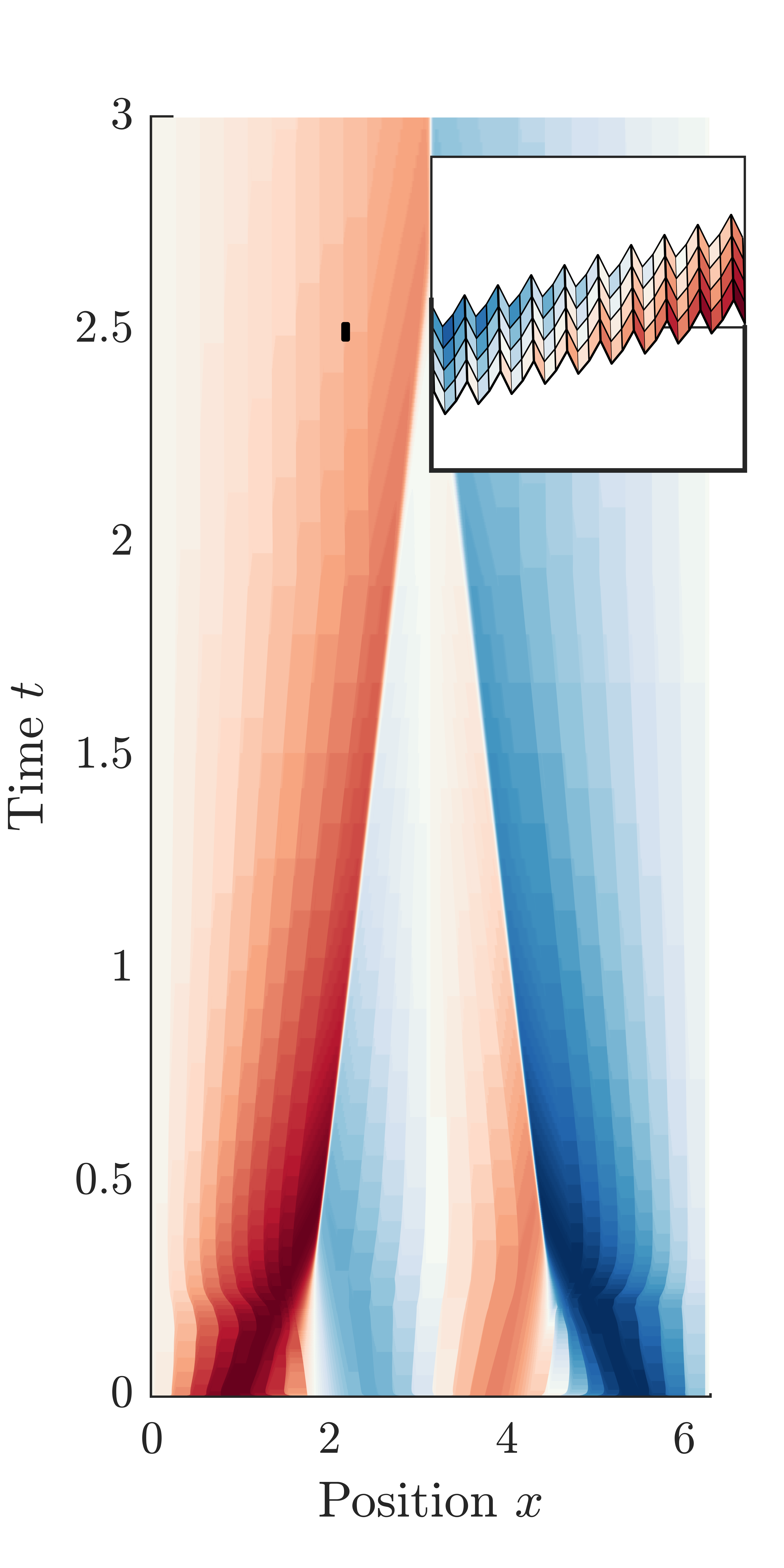}
& 
\includegraphics[scale=0.95]{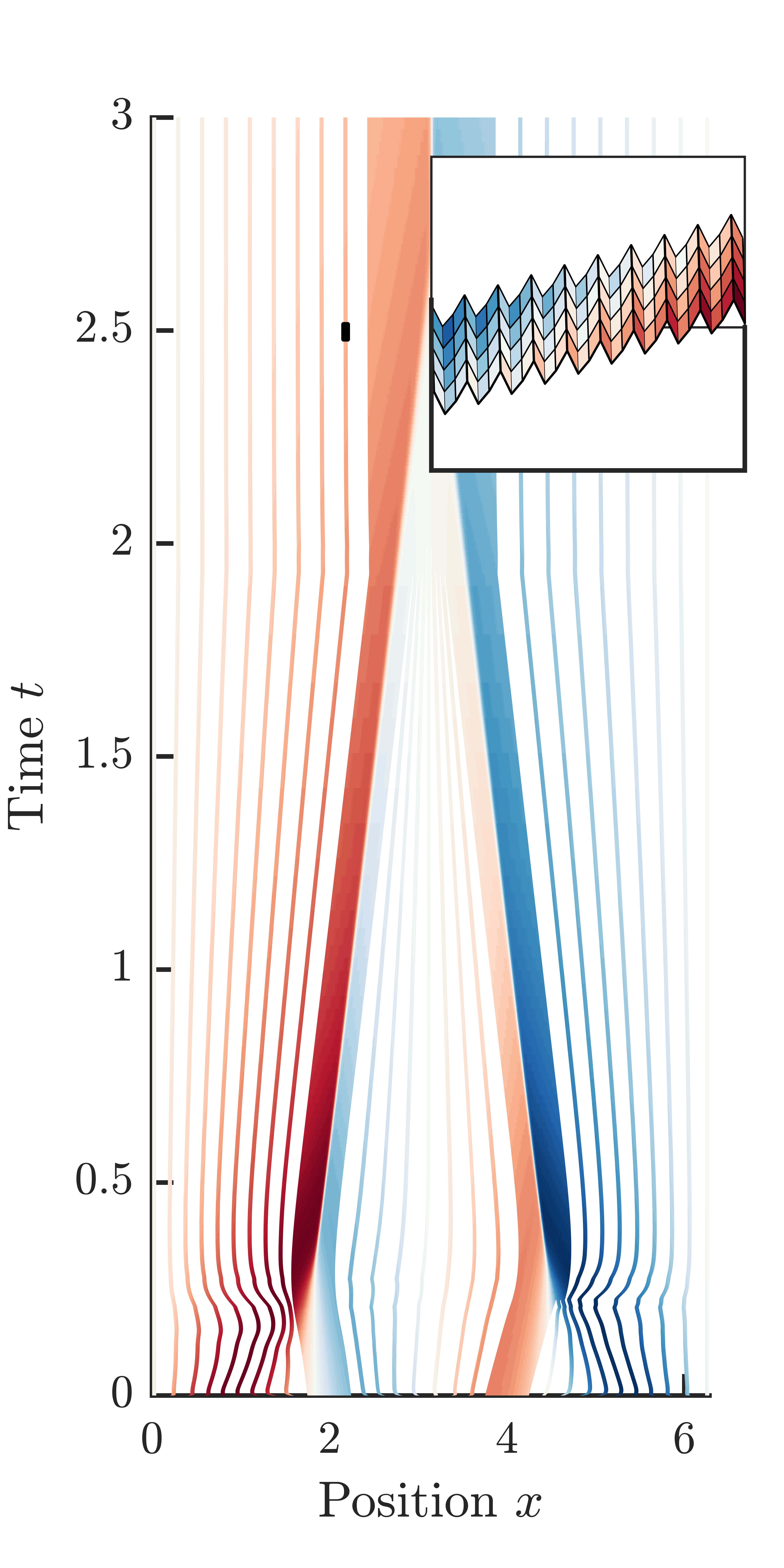}
\end{tabular}       
\caption{For the example of \cref{sec:e2}: \quad(a)~full domain simulation shows two shocks form at $t\approx0.3$ then merge together at $t\approx2.7$.
The insets show the detailed micro-scale structure near $(x,t)=(2,2.5)$.
\quad(b)~The patch simulation initially contains two \Pat{}es, but neither is aligned with the emerging shocks.
The \Pat\ initially located at $x=2$ moves to cover the left-hand shock (visible as a sudden change of colour from red to blue).
The \Pat\ initially at $x=4$ is further from the right-hand emerging shock: consequently the two nearest patches merge with the \Pat, empowering the wider \Pat\ to resolve the forming shock.  }
    \label{fig:hetDblSin}
\end{figure}

This example, shown in \cref{fig:hetDblSin}, illustrates that our merging and moving patch scheme moves \Pat{}es towards shocks when the \Pat{}es are initially misplaced.
Thereafter, \Pat{}es move with the shocks and merge (\cref{fig:mrg}(c)) with each other when necessary. 

We choose the spatial domain to be~\([0,2\pi]\), initial conditions $u(x,0) = \sin(2x)+\frac12\sin(x)$ sampled on the micro-grid with spacing $\dx=0.0016$, so that the full domain simulation of \cref{fig:hetDblSin}(a) has \(4\,000\)~micro-grid points in space, and choose boundary conditions $u=0$ at \(x=0,2\pi\).
From these initial conditions, two shocks form around time $t\approx0.3$, move towards each other, and eventually merge around $t\approx2.7$ (\cref{fig:hetDblSin}(a)).

\begin{table}
\caption{\label{tbltwo}specific heterogeneous coefficients, three-periodic, for the \ode{}s~\eqref{hetModel} used in the example of \protect\cref{sec:e2,fig:hetDblSin,fig:hetDblSinRmse}.}
\begin{equation*}
\begin{array}{r|lll}
\gamma_k& 3.14  &  0.37  &  0.39
\\
\epsilon_k& 0.0054  &  0.099  &  0.0096
\end{array}
\end{equation*}
\end{table}
We set the heterogeneous strengths $\sige=2$ and $\sigg=1$ with micro-grid periodicity $\kappa=3$, so that the micro-scale, with its log-normal coefficients, is strongly heterogeneous.
\cref{tbltwo} lists the specific realisation used here.
The listed and used diffusivities~\(\epsilon_k\) were scaled to have harmonic mean \(\epsilon=0.01\)\,.

\begin{figure}
\centering
\includegraphics[scale=0.95]{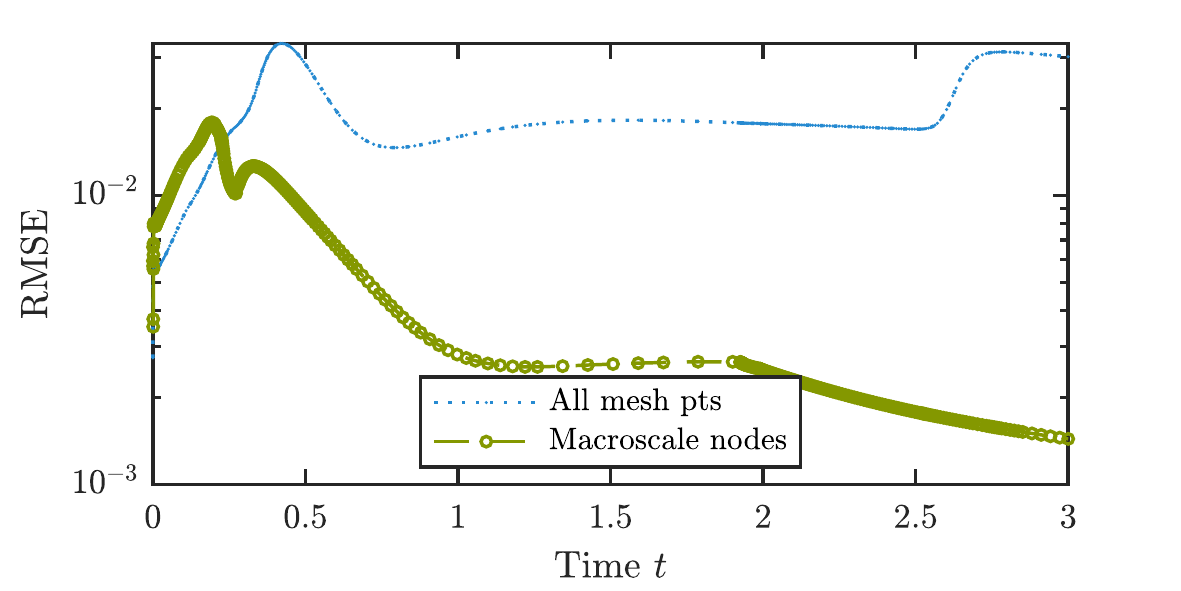} 
\caption{Errors (\textsc{rmse}) over time for the example of \cref{sec:e2}.
In this case shocks form around $t\approx0.3$ and merge $t\approx2.7$.
Both the micro-scale \textsc{rmse} (calculated at all patch micro-grid points) and the macro-scale \textsc{rmse} (calculated at the centre of each patch) remain controlled, but in particular the macro-scale \textsc{rmse} is excellent.
\cref{sec:cons} discusses consistency properties for moving patches in terms of the error at macro-scale nodes, and this low macro-scale \textsc{rmse} indicates our \text{scheme performs excellently.}}
    \label{fig:hetDblSinRmse}
\end{figure}

In order to stress-test the behaviour of the new \Pat{}es, we begin the simulation with two extant \Pat{}es---situated in the wrong locations to track the shocks.
Accurate simulation requires that the \Pat{}es move towards the forming~shocks, without destructively interfering with the simulated~solution.
We evenly space $28$~ordinary patches in the gaps between boundaries and \Pat{}es.
Each ordinary patch consists of $30$~micro-scale points, and each \Pat\ of~$300$.
All other parameters are as described in \cref{sec:e1}.
The simulation by patches is compared to an expensive, accurate, full-domain, simulation in \cref{fig:hetDblSin}.
The cross-section view (not shown) is similar to \cref{fig:sinhetDetails}: the micro-scale heterogeneity, and macro-scale structure, of the solution is well represented.
\cref{fig:hetDblSinRmse} confirms the patch schemes accuracy over all simulation times by plotting the Root Mean Squared Error (\textsc{rmse}) for both macro- and micro-scale variables.

\subsection{Example: travelling and colliding shocks}
\label{sec:e3}

\begin{figure}
\centering
\begin{tabular}{c c}
(a)&(b)\\
\includegraphics[scale=0.95]{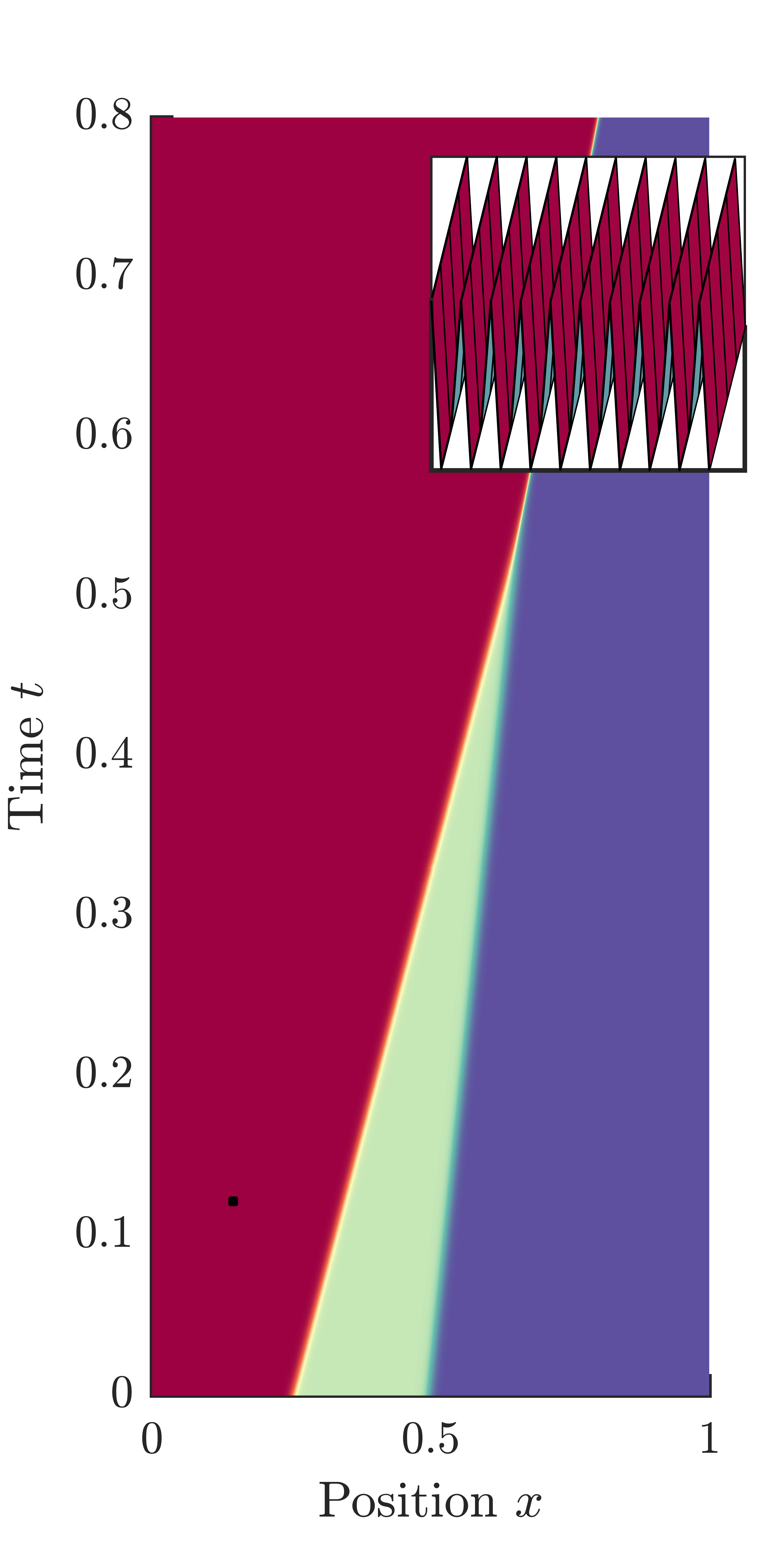}
& 
\includegraphics[scale=0.95]{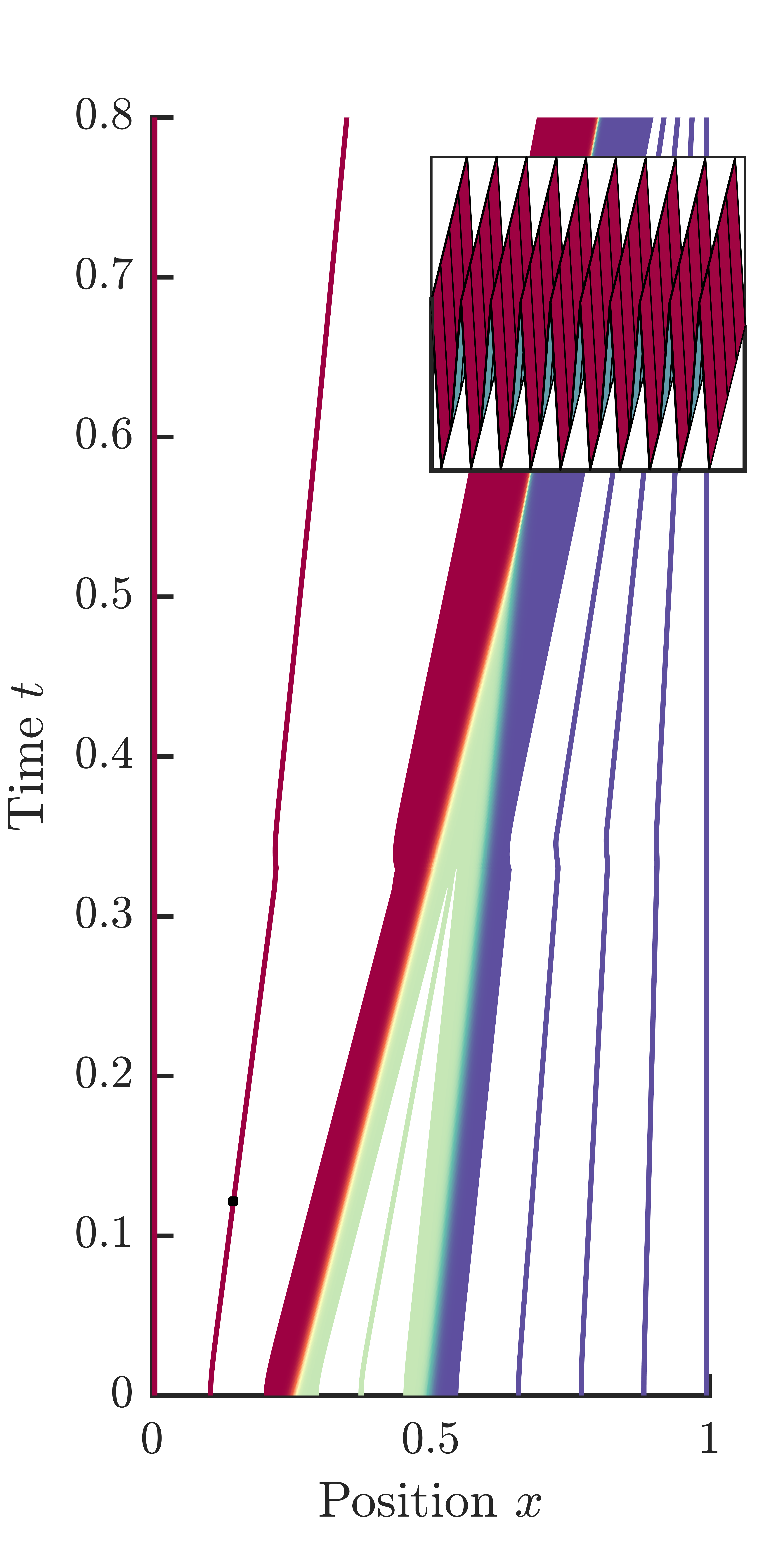}
\end{tabular}       
\caption{For the example of \cref{sec:e3}: \quad(a)~Overview of the full-domain simulation.
The main plot shows a rolling average of~$u$: regions coloured red (left side) are $u\approx1$, green (middle) are $u\approx0.5$, and dark blue (right side) are $u\approx0.1$.
\quad(b) Accurate patch simulation, initially consisting of a \Pat\ centred at $x=0.25$, a \Pat\ centred at $x=0.5$, and seven other patches. 
Insets detail the micro-scale structure near $(x,t)=(0.1,0.1)$, indicated with a black box on each.}
    \label{fig:hetDblShock}
\end{figure}

This example, shown in \cref{fig:hetDblShock}, simulates moving shocks between regions where the (averaged) macro-scale solution is almost constant in space, showing that the moving patch scheme does not destructively interfere with such piecewise \text{`constant' solutions.}

We adapt an example from \cite{Huang10}.
\emph{For homogeneous advection and diffusion} an exact solution to Burgers' \pde\ is
\begin{equation*}
u(x,t) = \frac 
{ 0.1\exp\left[\frac{0.5-x-4.95t}{20\epsilon}\right] 
+ 0.5\exp\left[\frac{0.5-x-0.75t}{4\epsilon}\right] 
+ \exp\left[\frac{0.375-x}{2\epsilon}\right] }
{ \exp\left[\frac{0.5-x-4.95t}{20\epsilon}\right] 
+ \exp\left[\frac{0.5-x-0.75t}{4\epsilon}\right] 
+ \exp\left[\frac{0.375-x}{2\epsilon}\right] } .
\end{equation*}
We provide initial and boundary conditions from this exact~solution to the heterogeneous Burgers'-like system~\eqref{hetModel}.
The solution consists of two shocks that move to the right at different speeds: the faster moving shock overtakes the slower at $t\approx0.55$.
The particular difficulty in simulation is that the macro-scale solution contains large regions of (almost) constant~$u$ (\cref{fig:hetDblShock}(a)).
Accurate simulation with moving patches (\cref{fig:hetDblShock}(b)) requires that the predicted field between the two fronts remain approximately constant.
Patch movements may distort the solution~$u$ via the chain rule in \cref{chainU}.
Together with the heterogeneous structure on each patch, a poorly implemented scheme would bend or warp the solution~$u$ \text{between the fronts.}

\begin{table}
\caption{\label{tblthree}specific heterogeneous coefficients, three-periodic, for the \ode{}s~\eqref{hetModel} used in the example of \protect\cref{sec:e3,fig:hetDblShock,fig:DblDetails,fig:Dblrmse}.}
\begin{equation*}
\begin{array}{r|lll}
\gamma_k& 1.0039  &  0.9948  &  1.0013
\\
\epsilon_k&0.0013  &  0.0005  &  0.019
\end{array}
\end{equation*}
\end{table}%
We choose the spatial domain to be~\([0,1]\), with a micro-grid of spacing $\dx=1/3000$, so that the full domain simulation of \cref{fig:hetDblShock}(a) has \(3\,000\)~micro-grid points in space.
We implement large heterogeneous diffusion $\sige=1$ and relatively small heterogeneous advection $\sigg=0.2$, with periodicity \text{$\kappa=3$}.
\cref{tblthree} lists the specific realisation used here.
The listed diffusivities~\(\epsilon_k\) were scaled to have harmonic mean \(\epsilon=0.001\) (one-tenth the harmonic mean of \text{the previous examples).}

\begin{figure}
\includegraphics{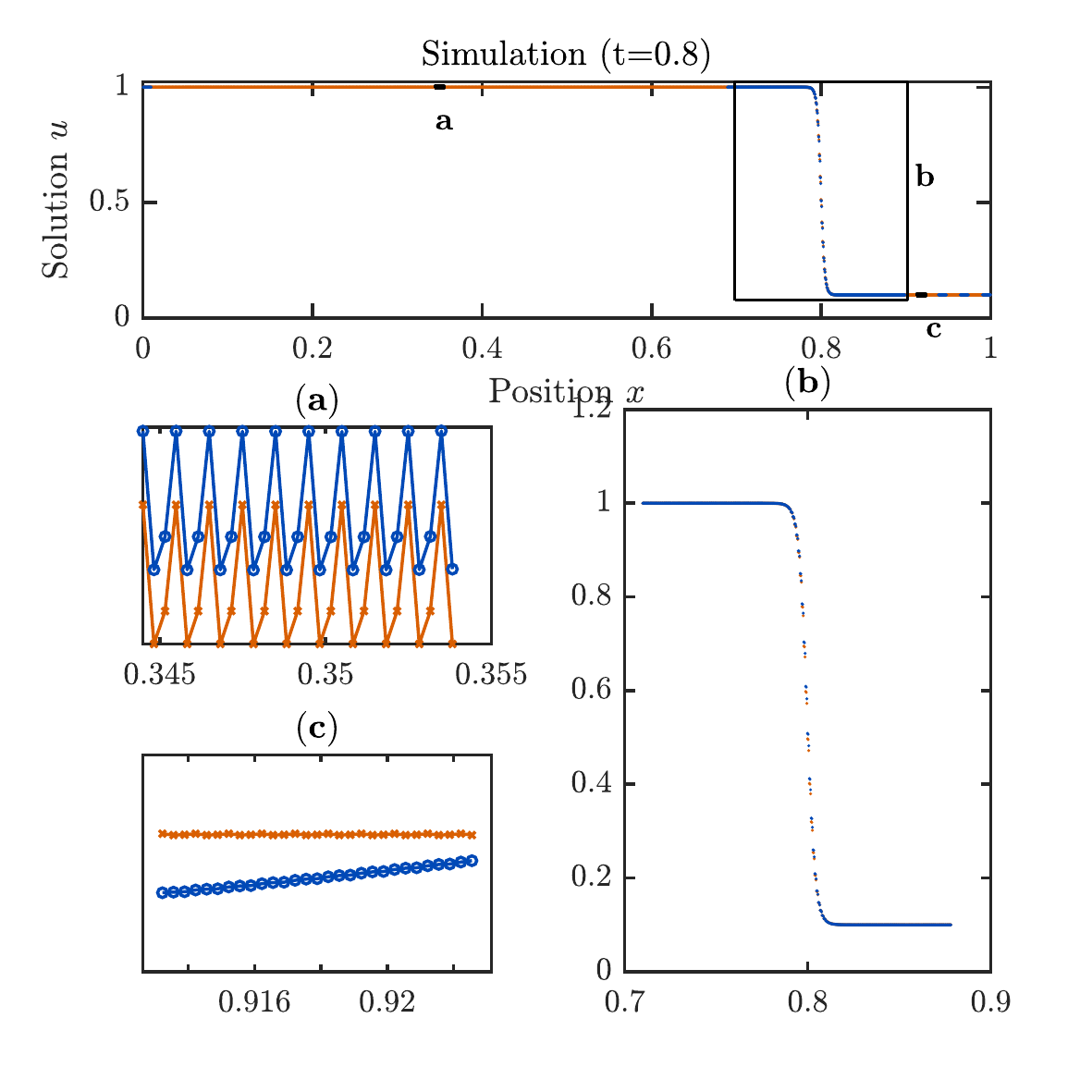}
\caption{Detailed view at the final time \(t=0.8\) of the example of \cref{sec:e3}.
There are small systematic errors on both sides of the shock. (a)~On the left side of the shock, the patch prediction~(blue) is slightly (about~$0.001$) above the full simulation~(orange). (c)~On the right side of the shock, the patch prediction is again slightly below the trusted simulation. (b)~the patch prediction of the shock structure is accurate on both the micro- and macro-scale. }
\label{fig:DblDetails}
\end{figure}

\begin{figure}
\includegraphics{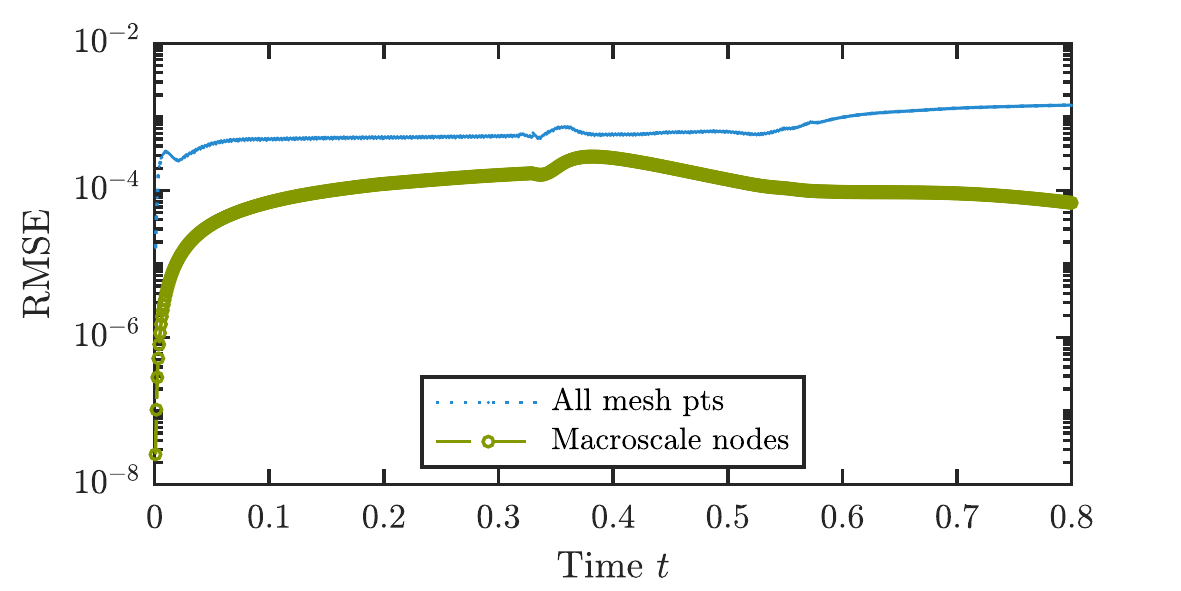}
\caption{Errors (\textsc{rmse}) over time of the patch scheme for the example of \cref{sec:e3}.
Both the micro-scale, and particularly the macro-scale, are accurately predicted by our patch scheme.  }
\label{fig:Dblrmse}
\end{figure}

We initialise simulation with two \Pat{}es centred over the initial shock locations, at $x=0.25$ and $x=0.5$.
Seven patches are placed evenly across the remaining space: each patch consists of $30$~micro-scale points and each \Pat\ has $300$~points.
As the shocks move rapidly across space~$x$, we set the mesh redistribution parameters to $\tau=1$ and $\beta=0.01$.
These rapid redistributions are precisely the situation that may induce errors in the solution~$u$ via \cref{chainU}.
To ameliorate such effects we couple the patches together with quadratic interpolation, $\Gamma=1$.
\Cref{fig:hetDblShock}(b) shows an excellent macro-scale prediction by our patch scheme.
\cref{fig:DblDetails} shows that the patches on either side of the shock have some small errors induced by patch movements, but these errors are only visible on the micro-scale.
The \textsc{rmse} in \cref{fig:Dblrmse} confirms that the patch predictions are accurate over all times.

\section{Emergence and consistency for moving and merging patches}
\label{sec:cons}

The three examples presented indicate that our moving patch scheme makes accurate predictions.
This section outlines some theoretical support for moving patches (\cref{sec:mov}) in terms of the emergence and consistency already proved for stationary patches.
Complete emergence and consistency proofs are deferred \text{for future research.}

\paragraph{Emergence of macro-scale dynamics for moving patches}
The patch scheme is a multiscale dynamical system:
the underlying system, such as~\eqref{hetModel}, has micro-scale lattice interactions on the length scale~\(\dx\);
the imposed patches have a length scale~\(2h\) where~\(h\) is a typical half-width \(h_j=n_j\dx\)\,;
the separation of patches has typical size~\(H\) representing the separations~\(H_j\);
and these are all placed into a spatial domain of some size~\(L\).
The dynamics are correspondingly multiscaled.
In the general class of systems where diffusion is the main dissipative mechanism, such as~\eqref{hetModel}, then the sub-patch dynamics typically dissipate in time at a rate faster than~\Ord{1/h^2}.
Conversely, the inter-patch dynamics dissipate at rates slower than~\Ord{1/H^2}.
Thus, upon neglecting fast transients, the emergent long-lasting dynamics of the patch scheme should correspond to the relatively slow inter-patch interactions.
That is, we expect the emergent slow manifold of the patch scheme to be that of the macro-scale inter-patch interactions. 

For a wide range of systems it has been proven that patch schemes with \(N\)~patches generally possess an \(N\)-dimensional emergent slow manifold, a manifold parametrised by one `macro-scale' variable per patch.
These proofs apply to systems in multi-D space with either homogeneous micro-scale \cite[e.g.]{Roberts2011a, Alotaibi2017a, Bunder2019d} or with heterogeneous micro-scale \cite[e.g.]{Bunder2013b, Bunder2020a}.
However, these proofs are predicated on the patches being stationary with identical spacing, which does not immediately hold for the moving \text{patches of \cref{sec:mov}.}

However, we now recast the moving patch system of \cref{sec:mov} to the requisite form for this theory to apply.
For simplicity, let's just address a micro-scale of the continuum \pde\ \(u_t=f(u,x)\) for a field~\(u(x,t)\), such as the equivalent \pde\ of~\eqref{hetModel}.
For moving patches we begin to recast the problem, as in \cref{sec:stdMov}, by seeking a field~\(\hat u(\xi,t)\) in terms of~\(\xi\), via some transform of space \(x=x(\xi,t)\).
The field~\(\hat u\) then satisfies the \pde~\eqref{chainU}, namely \(\hat u_t=f(\hat u,x)+(\hat u_\xi/x_\xi)x_t\)\, and the field~\(x\) satisfies the moving mesh \pde~\eqref{mmpde}.
Then we recast the \pde{}s for~\(\hat u\) and~\(x\) as one \pde\ system for the \(\RR^2\)~valued field $\uv(\xi,t):=\big(\hat u(\xi,t)\C x(\xi,t)\big)$.
For this \uv-system, the patches are stationary in~\(\xi\) with identical spacing.
Since the moving mesh \pde~\eqref{mmpde} is diffusive, then for suitable dissipative \pde{}s for~\(u\), the established emergence and consistency theories apply to the system for~\uv, and hence to the moving patch scheme \cite[e.g.]{Roberts2011a, Alotaibi2017a, Bunder2019d, Bunder2013b, Bunder2020a}.

With moving patches, the system possesses two fields, \(\hat u(\xi,t)\) and~\(x(\xi,t)\), evolving in time in \(\xi\)-space.
Consequently, we are generally assured that the moving patch scheme with \(N\)~patches, possesses emergent macro-scale dynamics that form the \(2N\)-dimensional slow manifold of the moving patch system.
The \(2N\)-D comes from \(N\)-D macro-scale dynamics of the field~\(u\), and \(N\)-D of the macro-scale dynamics the spatial transform~\(x(\xi,t)\) of the evolving \text{patch distribution.}

\paragraph{Macro-scale consistency for moving patches}
Similarly, the consistency of previous patch schemes has been established in diverse contexts: \cite{Roberts06d, Roberts2011a} proved consistency for patches in multi-D space, and \cite{Bunder2013b, Bunder2020a} proved consistency of two patch schemes for heterogeneous micro-scale structures.
These consistency results establish that a patch scheme with polynomial coupling~\eqref{patCoup} of order~\(2\Gamma\) generally has emergent macro-scale dynamics that are the same as the underlying micro-scale equations to errors~\Ord{H^{2\Gamma}}.
Hence one can control errors in the patch scheme by suitable choice of patch spacing~\(H\) and order~\(\Gamma\).

The arguments in the previous paragraphs assure us that for suitable dissipative \pde{}s for~\(u\), these established consistency results also apply to the moving patch scheme. 
The established consistency is in terms of the patch spacing~\(\hat H\) in~\(\xi\).
However, in most scenarios the patch spacing~\(H_j\) in~\(x\) is proportional to~\(\hat H\), and so the order of consistency is maintained.
Moreover, the moving patch scheme is designed to flexibly resolve regions of large curvature in~\(x\), so the coefficients of~\(H_j^{2\Gamma}\) in local error terms should, by design, make errors more commensurate across the spatial domain.
Consequently, we reasonably conjecture that the~\Ord{H^{2\Gamma}} consistency results established for patch schemes also apply to the \text{moving patch scheme of \cref{sec:mov}.}

The above two appeals to existing theories are clearly not complete as most details are missing, and we have not yet established that the proofs extend to patches with differing sizes.
Nonetheless, the \emph{outline} presented here supports our conjecture that the moving patch scheme has emergent macro-scale dynamics which is consistent with that of the underlying \text{micro-scale physics.}

\paragraph{Emergence and consistency for \emph{merging} patches}
The above discussion addresses moving patches, but not the merging of patches, nor simulations with \Pat{}es.
First, address \Pat{}es, and suppose there are \(M\)~\Pat{}es.  
These effectively divide the spatial domain into \((M+1)\)~regions, each region with the moving patch scheme, and so we expect each region to have appropriate emergent consistent macro-scale dynamics. 
The regions are smoothly separated by the interior of a \Pat\ in which the full micro-scale details, such as~\eqref{hetModel}, are computed, and in which the transition layers between the regions are fully resolved.
Consequently, we expect the simulation over the entire domain to have emergent macro-scale dynamics that consistent with the \text{micro-scale physics.}

Second, consider the merger of patches (\cref{sec:patMrg}), say merging patches~\(s\) and~\(s+1\) at time~\(t'\).
We argue that our scheme is nearly continuous and so there is a smooth enough transition from before to after the merge.
For all micro-grid points not adjacent to a patch edge, time derivatives are \emph{identical} at times~$t'-$ and~$t'+$ as patch merges do not alter the location of or field at any patch-interior micro-grid points.
The only issue is how the merge affects edge-patch values. 
For patches~\(j\) with \(|j-s-\tfrac12|\geq\Gamma+\tfrac12\) the interpolation~\eqref{patCoup} is unaffected by the merge, and so edge-patch values are identical through the merge.
For nearer patches, the interpolation changes from one scheme that has errors~\ord{H^{2\Gamma}}, measured over the (sub-domain) regions, to another scheme that also has global errors~\ord{H^{2\Gamma}} \cite[e.g.]{Shoosmith75, Beyn79}.
Hence the edge-patch values only change by an amount~\ord{H^{2\Gamma}} which is small enough to ensure changes in next-to-edge \(\de tu_{j,\pm(n_j-1)}\) across the merge are smooth enough to preserve the order of error in the scheme.
Similarly, at the point of collision, the right-edge of the \(s\)th~patch and the left-edge of the \((s+1)\)th~patch have \(u\)-values with errors~\ord{H^{2\Gamma}}, and so their average used in the merge also has error~\ord{H^{2\Gamma}}.
Consequently, the proposed merge is smooth enough to preserve the order of error, and hence the consistency with the underlying given micro-scale system, in \text{the moving patch scheme.}

\section{Discussion}
\label{sec:disc}

We have empowered the Equation-Free patch scheme to detect and track macro-scale features in 1D space, for example shocks, fronts, or cracks.
Our modification to the patch scheme comes from adaptive moving mesh schemes, combined with a strategy to merge patches on-the-fly.
This modification has good theoretical support: moving patches are justified by considering the moving mesh as another set of dependent variables, simulated on a fixed and uniform computational mesh on the index-space.
The moving patches therefore correspond to stationary patches on the computational mesh, and the typical consistency arguments for patch schemes follow.
We justify the merging of patches by demonstrating that patch merges preserve properties of the solution and of the governing dynamical equations: in particular, the phase structure of periodic heterogeneous micro-scales is preserved across merges.

The novel patch scheme was thoroughly tested on a modified 1D Burgers' equation with heterogeneous advection and heterogeneous diffusion terms.
We demonstrate accurate simulation of several difficult problems.
In particular, our scheme outperforms neural networks (even when the neural nets are only solving problems with a homogeneous micro-scale). 

In future work we plan to augment this scheme with a strategy to allow \Pat{}es to separate, increasing the flexibility of the scheme for situations where shocks decay and evaporate.
For example, when the inter-patch interpolants from the two sides of a \Pat\ agree, to some tolerance, in the interior of the \Pat\ then it may be split into two ordinary patches.

\paragraph{Towards multiple spatial dimensions} 
In this article the scheme and applications have been restricted to 1D~space.
We plan further research and development to extend the scheme to more spatial dimensions. 
The foundational moving patch scheme in \cref{sec:movPat} extends to multi-D without difficulty, using the appropriate formulation of both the patch scheme and established moving mesh methods. 
The outstanding challenge is to extend the details of a \Pat\ and the merging of patches in multiple dimensions, while retaining both accuracy and efficiency.

\paragraph{Acknowledgements}
This research was funded by the Australian Research Council under grants DP200103097 and DP180100050.  
The work of I.G.K. was also partially supported by the DARPA PAI program.

\end{document}